\documentclass[reqno,a4paper]{amsart}

\usepackage{url}
\usepackage{amsmath,amsthm,amssymb}
\usepackage{amscd}
\usepackage{amsfonts}
\usepackage{xcolor}
\usepackage[german,french,english]{babel}
\usepackage{setspace}
\makeatletter
 
  \@addtoreset{equation}{section}
\makeatother
\newtheorem{thm}{Theorem}[section]
\newtheorem{dfn}[thm]{Definition}
\newtheorem{lem}[thm]{Lemma}

\newtheorem{conj}[thm]{Conjecture}

\makeatletter
\@namedef{subjclassname@2020}{%
  \textup{2020} Mathematics Subject Classification}
\makeatother

\begin{document}

%\doublespacing
%%%%%%%%%%%%%%%%%%%%%%%%%%%%%%%%%%%%%%%%%%%%%%%%%%%%%%%%%%%
\title{On the Non-extendibility of Triangular $D(n)$-Pairs $\{n,2n\}$ with $n = 2k^2 \pm k$}

\author{Alan Filipin}
\address{Faculty of Civil Engineering, University of Zagreb, Fra Andrije Ka\v{c}i\'{c}a-Mio\v{s}i\'{c}a 26, 10000 Zagreb, Croatia}
\email{filipin@grad.hr}

\author{Yasutsugu Fujita}
\address{Department of Mathematics, College of Industrial Technology, Nihon University, 2-11-1 Shin-ei, Narashino, Chiba, Japan}
\email{fujita.yasutsugu@nihon-u.ac.jp}

\author{Ana Jurasi\'{c}}
\address{Faculty of Mathematics, University of Rijeka, Radmile Matej\v{c}i\'{c} 2, 51000 Rijeka, Croatia}
\email{ajurasic@math.uniri.hr}

\date{\today}

%%%%%%%%%%%%%%%%%%%%%%%%%%%%%%%%%%%%%%%%%%%%%%%%%%%%%%%%%%%
\begin{abstract}

We study a class of Diophantine problems arising from triangular analogues of classical Diophantine tuples. A set of positive integers is called a triangular $D(n)$-tuple if the product of any two distinct elements, increased by $n$, is a triangular number. In this paper, we investigate the extensibility of triangular $D(n)$-pairs of the form $\{n,2n\}$, where $n = 2k^2 + \varepsilon k$, with $k$ an odd positive integer and $\varepsilon \in \{-1,1\}$ such that $4k + \varepsilon$ is prime. We prove that for all $(k,\varepsilon)$ in a finite explicitly determined set, the pair $\{n,2n\}$ cannot be extended to a triangular $D(n)$-quadruple. The obtained results give further evidence for a general non-extensibility phenomenon, which we formulate as a conjecture.
\end{abstract}


\maketitle

%%%%%%%%%%%%%%%%%%%%%%%%%%%%%%%%%%%%%%%%%%%%%%%%%%%%%%%%%%%
\section{Introduction}\label{sec:int}
%%%%%%%%%%%%%%%%%%%%%%%%%%%%%%%%%%%%%%%%%%%%%%%%%%%%%%%%%%%
%%%

%%%%%
The study of Diophantine tuples originates in the classical work of Diophantus, who exhibited sets of rational numbers with the property that the product of any two distinct elements, increased by $1$, is a perfect square. This problem has since inspired numerous generalizations, particularly in the direction of replacing squares by other figurate numbers. One natural generalization considers triangular numbers $\Delta_k = \frac{k(k+1)}{2}$, $k\in\mathbb{N}$, (sequence A000217 in the OEIS \cite{oeis}).
\begin{dfn}\label{df1}
Let $n$ be an integer. A set $\{a_{1}, \ldots ,a_{m}\}$ of $m$ distinct positive integers such that
$a_{i}a_{j}+n$ is a triangular number, for $1\leq i<j\leq m$, is called a triangular Diophantine $m$-tuple with the property $D(n)$ or simply a Triangular $D(n)$-$m$-tuple. \end{dfn}

In analogy with the classical case, a central problem is to determine the maximal possible size of such sets and to understand the conditions under which a given tuple can be extended.  Juki\'{c} Bokun and Soldo \cite{JBiS}  considered Triangular $D(-1)$-tuples and proved that the only Triangular $D(-1)$-triples of the form $\{1, 2, c\}$, $c = 2^n p$, where $n$ is a non-negative integer and $p$ is a prime, are those with $c \in \{11, 46, 352, 11936\}$. They proved that for these $c$’s no Triangular $D(-1)$-quadruple of the form $\{1, 2, c, d\}$ exists.  Several authors have shown that specific pairs of the form $\{n,2n\}$ cannot be extended to quadruples for certain values of $n$. These results rely on transforming the problem into systems of Pellian equations and, in more difficult cases, to determining integral points on elliptic curves. Filipin and Szalay \cite{FS} studied the extensibility of the Triangular $D(1)$-pair $\{1,2\}$, and proved that it cannot be extended to a Triangular $D(1)$-quadruple $\{1,2,c,d\}$. Filipin, Jurasi\'{c}, and Szalay \cite{FJS} proved that the Triangular $D(3)$-pair $\{3,6\}$ cannot be extended to a Triangular $D(3)$-quadruple $\{3,6,c,d\}$. Also, Filipin, Jurasi\'{c}, and Togb\' e \cite{FJT} proved that the Triangular $D(n)$-pair $\{n, 2n\}$, for $n=10$ and $n=28$, cannot be extended to a quadruple. The problem described in \cite{FJT} doesn't belong to the generalization we are dealing with, because $10=2\cdot2^2+2$ and $28=2\cdot4^2-4$, but $4\cdot2+1=9$ and $4\cdot4-1=15$ are not primes and $k$ is even. Motivation for this paper was to generalize such results and to study the extension of Triangular $D(n)$-pair of the form $\{n, 2n\}$. However, such problem is not straight-forward to solve. The main problem is finding fundamental solutions of Pellian equations obtained by transforming our problem to solving the system of such equations. More precisely, the main problem was searching of integral points on elliptic curves associated to the problem as we will explain later. Thus, we have to restrict ourselves to the set of the pairs $(k,\varepsilon)$, suitable for solving, for now. We consider the family
\[
n = 2k^2 + \varepsilon k,
\]
where $k$ is an odd positive integer and $\varepsilon \in \{-1,1\}$, under the additional condition that $4k + \varepsilon$ is prime. This condition ensures a favorable arithmetic structure that allows us to treat several cases simultaneously. Our goal is to investigate whether the corresponding triangular $D(n)$-pair $\{n,2n\}$ can be extended to a quadruple. The main result of our work is the following theorem.

\begin{thm}\label{thm:main}
Let $k$ be an odd positive integer with $4k+\varepsilon$ prime, where $\varepsilon \in \{-1,1\}$. There does not exist a Triangular $D(2k^2+\varepsilon k)$-quadruple of the form $\{2k^2+\varepsilon k,2(2k^2+\varepsilon k),c,d\}$, where
\begin{eqnarray}(k,\varepsilon)\in S:=\{(3, -1), (3, 1), (5, -1), (7, 1), (9, 1), (11, -1), (13, 1), (27,1) \}.\label{thm:main-S}
\end{eqnarray}
\end{thm}

The proof is based on reducing the problem to a system of generalized Pellian equations and analyzing their solutions via recurrence sequences. A key step is the connection between the parameters involved and generalized pentagonal numbers, which provides a structural constraint on possible extensions. In addition, the problem leads naturally to the study of integral points on certain elliptic curves, whose analysis imposes further restrictions. Although our result is presently limited to a finite set of parameters, the method developed here applies uniformly across the family and highlights the main obstacles to a complete solution. This leads us to formulate the following conjecture: \begin{conj}\label{conjecture}
Let $k$ be an odd positive integer with $4k+\varepsilon$ prime, where $\varepsilon \in \{-1,1\}$.
Then, there does not exist a Triangular $D(2k^2+\varepsilon k)$-quadruple of the form $\{2k^2+\varepsilon k,2(2k^2+\varepsilon k),c,d\}$.
\end{conj}

\bigskip

%%%%%
This is a generalization of the cases where $k=1$, observed in \cite{FS} for $\varepsilon =-1$ and in \cite{FJS} for  $\varepsilon =1$. We, therefore, may assume that $k\geq 3$ and $4k+\varepsilon\geq 11$. Since $4k+\varepsilon$ is a prime, we know that \begin{equation}\label{prime}4k+\varepsilon=6p\pm 1,\end{equation} for some $p\in\mathbb{N}$.

\textbf{Case a)} If $4k=6p$, then $k$ is divisible by $3$, i.e., $k=3k_1$ for some  odd $k_1\in\mathbb{N}$ ($p=2k_1$)  and we are considering  a Triangular $D(3k_1(6k_1+\varepsilon))$-quadruple of the form $$\{3k_1(6k_1+\varepsilon) ,6k_1(6k_1+\varepsilon),c,d\}.$$ In this case $$4k+\varepsilon=12k_1+ \varepsilon,$$ for some  odd $k_1\in\mathbb{N}$. More precisely, $k=6k_2+3$ and $4k+\varepsilon$ is a prime number of the form $24k_2+11$ or $24k_2+13$, where $k_2$ is a non-negative integer.

\textbf{Case b)} If $4k=6p- 2\varepsilon$, then $2k+\varepsilon=3p$ so $k=3k_1+\varepsilon$, $p=2k_1+\varepsilon$ for some  even $k_1\in\mathbb{N}$ and $\gcd (p,k)=1$. We are considering a Triangular $D(3(2k_1+\varepsilon)(3k_1+\varepsilon))$-quadruple of the form \begin{equation}\{3(2k_1+\varepsilon)(3k_1+\varepsilon),6(2k_1+\varepsilon)(3k_1+\varepsilon),c,d\}.\end{equation} In this case $$4k+\varepsilon=12k_1+ 5\varepsilon,$$ for some  even $k_1\in\mathbb{N}$. For $k_1=0$ and $\varepsilon=1$, we would have the case observed in \cite{FJS}. More precisely, $k=6k_2\pm 1$ and $4k+\varepsilon$ is a prime number of the form $24k_2\pm 5$, where $k_2\in\mathbb{N}$.

We can observe both cases together, by taking \begin{equation}\label{lambda}2k^2+\varepsilon k=3\lambda,\end{equation}where $$\lambda=k_1(6k_1+\varepsilon),$$ for some  odd $k_1\in\mathbb{N}$, in the case a) and $$\lambda=(2k_1+\varepsilon)(3k_1+\varepsilon),$$ for some even $k_1\in\mathbb{N}$, in the case b). From (\ref{lambda}), we further obtain \begin{equation}\label{square}1+24\lambda=(4k+\varepsilon)^2.\end{equation}  It is easy to see that $\lambda$ is an odd integer. Note that if $\varepsilon =1$,  in the case b) for $\lambda$ we obtain every other member of the sequence of pentagonal numbers with odd indices (subsequence of the sequence A033570 in the OEIS \cite{oeis}).  Note also that the lowest possible value of $\lambda$ is $\lambda=5$, obtained in the case a) for $k_1=1$ and $\varepsilon=-1$. Moreover, from (\ref{square}), using (\ref{prime}) we obtain that \begin{equation}\label{lambda-pent}\lambda=\frac{p(3p\pm 1)}{2},\end{equation}with the sign $+$ in the case a) and the sign $-$ in the case b), so $\lambda$ is a generalized pentagonal number. Also, $3\lambda$ is a triangular number.

For the proof of the Theorem \ref{thm:main}, we start by following the strategy used in \cite{FS} and \cite{FJS}. We observe the general form of the problem and, firstly, we find the sequence $(c_l)_{l\in\mathbb{N}}$ of possible values of $c$ such that  $\{2k^2+\varepsilon k, 2(2k^2+\varepsilon k),c\}$ is a Triangular $D(2k^2+\varepsilon k)$-triple. In Section~2 we derive the associated system of generalized Pellian equations and describe the structure of their solutions. Section~3 establishes auxiliary properties of the parameters and their relation to figurate numbers. In Section~4 we analyze the fundamental solutions. In Sections~5--7 we complete the proof of the main theorem by combining congruence arguments, properties of the associated recurrence sequences, and estimates arising from linear forms in logarithms to exclude all possible extensions.

%%%%%%%%%%%%%%%%%%%%%%%%%%%%%%%%%%%%%%%%%%%%%%%%%%%%%%%%%%%
\section{System of generalised Pellian equations}
%%%%%%%%%%%%%%%%%%%%%%%%%%%%%%%%%%%%%%%%%%%%%%%%%%%%%%%%%%%
By considering a general form of the problem, we will be able to observe all elements of the set $S$ from Theorem \ref{thm:main} together and we will also  make theoretical basis useful for eventually future proving of Conjecture \ref{conjecture}.

Let $k$ be an odd positive integer and let $\varepsilon \in \{-1,1\}$. It is easy to see that  $\{2k^2+\varepsilon k,2(2k^2+\varepsilon k)\}$ is a Triangular $D(2k^2+\varepsilon k)$-pair.
Assume that $\{2k^2+\varepsilon k,2(2k^2+\varepsilon k),c\}$
%where $c>2(2k^2+\varepsilon k)$,
is a Triangular $D(2k^2+\varepsilon k)$-triple with the smallest possible $c\neq 0$. Similar as in classical problem of Diophantus, we have the following definition.
\begin{dfn}A Triangular $D(n)$-$m$-tuple in which elements are not necessarily distinct and may include zero is called improper Triangular $D(n)$-$m$-tuple.
\end{dfn}
Therefore, we assumed that we have a proper Triangular Diophantine triple and, in our case, there exist positive integers $s$ and $t$ such that
\begin{align*}
(2k^2+\varepsilon k)c+2k^2+\varepsilon k&=\frac{s(s+1)}{2},\\
2(2k^2+\varepsilon k)c+2k^2+\varepsilon k&=\frac{t(t+1)}{2},
\end{align*}
which yield
\begin{align}
8(2k^2+\varepsilon k)c+(4k+\varepsilon)^2&=S^2,\label{S}\\
16(2k^2+\varepsilon k)c+(4k+\varepsilon)^2&=T^2,\label{T}
\end{align}
where $S=2s+1$ and $T=2t+1$.
Eliminating $c$ from \eqref{S} and \eqref{T}, we obtain the Pellian equation
\begin{align}\label{TS}
T^2-2S^2=-(4k+\varepsilon)^2.
\end{align}
The solution $(T, S)$ is called primitive if $\gcd (T,S)=1$. If \eqref{TS} has a primitive integer solution $(T,S)$, then we easily see that
\begin{align*}
T=|A^2-B^2+2\mu AB|,\quad 4k+\varepsilon=|A^2-B^2-2\mu AB|,\quad S=A^2+B^2,
\end{align*}
where $A$, $B$ are coprime positive integers and $\mu \in \{-1,1\}$.
Hence, we would have $4k+\varepsilon \!\!\mod 8 \in \{1,2,6,7\}$, which contradicts the assumption that $k$ is odd.
Thus, \eqref{TS} has no primitive integer solutions $(T,S)$.

Suppose that $4k+ \varepsilon$ is a prime.
Then, we may write $T=(4k+\varepsilon)T'$ and $S=(4k+\varepsilon)S'$ for some positive integers $T'$ and $S'$, which leads to
\begin{align}\label{TS'}
(T')^2 - 2(S')^2 = -1.
\end{align}
Since it is clear that $\gcd(2(2k^2+\varepsilon k),4k+\varepsilon)=1$, by  \eqref{S} and \eqref{T} we also have
\begin{align}
8(2k^2+\varepsilon k)c'+1&=(S')^2,\label{S'}\\
16(2k^2+\varepsilon k)c'+1&=(T')^2,\label{T'}
\end{align}
where $c'=c/(4k+\varepsilon)^2$ is a positive integer.

Assume now that $\{2k^2+\varepsilon k,2(2k^2+\varepsilon k),c,d\}$, with $d>c$, is a Triangular $D(2k^2+\varepsilon k)$-quadruple.
Then, there exist odd integers $X$, $Y$ and $Z$ such that
\begin{align}
8(2k^2+\varepsilon k)d+(4k+\varepsilon)^2&=X^2,\label{X}\\
16(2k^2+\varepsilon k)d+(4k+\varepsilon)^2&=Y^2,\label{Y}\\
8cd+(4k+\varepsilon)^2&=Z^2.\label{Z}
\end{align}
In the same way as above, from \eqref{X}, \eqref{Y} and the fact that $4k+\varepsilon$ is a prime, we deduce that
\begin{align}
(Y')^2-2(X')^2&=-1,\label{YX'}\\
8(2k^2+\varepsilon k)d'+1&=(X')^2,\label{X'}\\
16(2k^2+\varepsilon k)d'+1&=(Y')^2,\label{Y'}
\end{align}
where $X'=X/(4k+\varepsilon)$, $Y'=Y/(4k+\varepsilon)$ and $d'=d/(4k+\varepsilon)^2$ are positive integers.
Hence, by \eqref{Z} we see that $Z'=Z/(4k+\varepsilon)$ is also a positive integer and
\begin{align}
8(4k+\varepsilon)^2c'd'+1=(Z')^2. \label{Z'}
\end{align}

%By using \cite[Theorem 108 a.]{nagel}, we get that $(Y',X')=(\pm 1,1)$ is the only fundamental solution of the equation (\ref{YX'}). Therefore, we have solutions $(Y',X')$, such that
It is well known that all solutions of Pellian equation \eqref{YX'} are given by
\begin{equation}\label{sedam'}Y'+X'\sqrt{2}=(1+\sqrt{2})(3+2\sqrt{2})^l,\ \ l\geq 0, \end{equation}
i.e.,
\begin{equation}\label{sedam}Y'+X'\sqrt{2}=(1+\sqrt{2})(T_l+U_l\sqrt{2}),\ \ l\geq 0,\end{equation}
where $(T_l,U_l)$ are defined by
\begin{eqnarray}\label{T_0}T_0&=&1,\ \ T_1=3,\ \ T_{l+2}=6T_{l+1}-T_l, \ \ l\geq 0,\\
\label{U}U_0&=&0,\ \ U_1=2,\ \ U_{l+2}=6U_{l+1}-U_l, \ \ l\geq 0.\end{eqnarray}
By \eqref{sedam}, we have
\begin{equation}\label{osam}(Y',X')=(T_l+2U_l, U_l+T_l).\end{equation}Note that $(U_l)_{l\geq 0}$ is a subsequence of the sequence of Pell numbers (sequence A000129 in the OEIS \cite{oeis}), i.e., its every other member. Therefore, we can use the formula for Pell numbers and express $$U_l=\frac{(1+\sqrt{2})^{2l}-(1-\sqrt{2})^{2l}}{2\sqrt{2}}\sim\frac{1}{2\sqrt{2}}(1+\sqrt{2})^{2l}, \ \ l\geq 0.$$ Also, $$U_l=P_l((1+\sqrt{2})^{l}+(1-\sqrt{2})^{l}),\ \ l\geq 0,$$ where $(P_l)_{l\geq 0}$ is a sequence of Pell numbers. For the sequence $(T_l)_{l\geq 0}$ (the sequence A001541  in the OEIS \cite{oeis}), it holds $$T_l=\frac{(3+2\sqrt{2})^{l}+(3-2\sqrt{2})^{l}}{2\sqrt{2}}\sim\frac{1}{2}(1+\sqrt{2})^{2l}, \ \ l\geq 0.$$

As possible values of $c'$ and $d'$, we get
\begin{equation}\label{devet}\frac{(U_l+T_l)^2-1}{8(2k^2+\varepsilon k)}=\frac{(U_l+T_l)^2-1}{24\lambda}.\end{equation}
By (\ref{T_0}) and (\ref{U}), we have \begin{equation}\label{sequence1}(U_l+T_l)_{l\geq 0}=1,5,29,169,985,...,\end{equation}  which is the sequence A001653 in OEIS \cite{oeis}.
% and $(-U_l+T_l)_{l\geq 0}=1,1,5,29,169,985,...$ so, as possible $c'$ and $d'$, it is enough to observe
Using (\ref{devet}) and (\ref{lambda}), we obtain
\begin{equation}\label{C_l}c'=c'_l=\frac{C_l}{\lambda}\end{equation}
and
\begin{equation}\label{fraction}C_l=\frac{(U_l+T_l)^2-1}{24}.\end{equation}
Therefore, \begin{equation}\label{sequence2}C_l=0,1,35,1190,40426, 1373295,...,\end{equation} which is the sequence A305539 in OEIS \cite{oeis}, i.e., the sequence of generalized pentagonal numbers which multiplied with $2$ are again generalized pentagonal numbers. Elements $d'=d'_l$ also have the form (\ref{C_l}). Moreover, using (\ref{square}), we conclude that the possible Triangular $D(3\lambda)$-triples have the form
\begin{equation}\label{deset}\Big\{3\lambda,6\lambda,24C_l+\frac{C_l}{\lambda}\Big\},\ \ \textup{for}\ \ l=0,1,2,....\end{equation}
Obviously, suitable candidates for $C_l$ are those which are multiples of $\lambda$, i.e., $C_l$ are generalized pentagonal numbers (sequence A001318 in the OEIS \cite{oeis}) and also multiples of an odd generalized pentagonal number $\lambda$ and we have
\begin{align*}c=(24\lambda+1)c'.\end{align*}
Similarly,
\begin{align*}d=(24\lambda+1)d'.\end{align*}

Using (\ref{lambda}), (\ref{square}) and (\ref{C_l}), we can transform equations (\ref{S'}) and (\ref{T'}) into
\begin{align}
24\lambda c'+1&=(S')^2,\label{S'_1}\\
48\lambda c'+1&=(T')^2,\label{T'_1}
\end{align}
and equations \eqref{X'}, \eqref{Y'} and \eqref{Z'} into
\begin{align}
24\lambda d'+1&=(X')^2,\label{X'_1}\\
48\lambda d'+1&=(Y')^2,\label{Y'_1}\\
8(1+24\lambda)c'd'+1&=(Z')^2. \label{Z'_1}
\end{align}

Eliminating $d'$ and using (\ref{square}) and the fact that  $X'=X/(4k+\varepsilon)$, $Y'=Y/(4k+\varepsilon)$  are integers, from equations \eqref{X'_1}, \eqref{Y'_1} and \eqref{Z'_1}, we obtain the system of generalized Pellian equations
\begin{align}
3\lambda(Z')^2- c'X^2&=3\lambda-c,\label{ZX'}\\
6\lambda(Z')^2-c'Y^2&=6\lambda-c.\label{ZY'}
\end{align}

\section{Additional properties of the parameters involved in the problem}

Let us denote $\lambda c'=C_l:=C$. By (\ref{S'_1}) and (\ref{T'_1}), \begin{align}
24C+1&=(S')^2,\label{S'_11}\\
48C+1&=(T')^2.\label{T'_11}\end{align} %Since $\lambda\geq 5$, from (\ref{C_l}), it follows \begin{equation}\label{c'-bound}c'\leq 0.2C\end{equation} and, from (\ref{deset}), it follows \begin{equation}\label{c-bound}c\leq 24.2C.\end{equation}

Knowing that $C$ is generalized pentagonal number (and, therefore, $3C$ is a triangular number), we have \begin{equation}\label{C-pent}C=\frac{m(3m\pm 1)}{2},\end{equation}from which we get $$24C+1=(6m\pm 1)^2,$$for some $m\in\mathbb{N}$. Therefore, by (\ref{C_l}), (\ref{S'_1}) and (\ref{X'_1}), \begin{equation}\label{s-mod}S'\equiv \pm 1  {\ ({\rm mod}\ {6})}\end{equation} and$$X'\equiv \pm 1  {\ ({\rm mod}\ {6})}.$$Since $2C$ is also pentagonal number, similarly, by  (\ref{C_l}), (\ref{T'_1}) and (\ref{Y'_1}), \begin{equation}\label{t-mod}T'\equiv \pm 1  {\ ({\rm mod}\ {6})}\end{equation} and$$Y'\equiv \pm 1  {\ ({\rm mod}\ {6})}.$$
More precisely, $S'$ and $X'$ are elements of the sequence $(U_l+T_l)_{l\geq 0}$ and, by (\ref{T_0})-(\ref{U}),
\[
U_l+T_l\equiv \begin{cases}1 \pmod{6}&\text{if~\,}l \equiv 0,3\pmod{4},\\
-1\pmod{6}&\text{if~\,}l \equiv 1,2\pmod{4}.\end{cases}
\]
%$U_0+T_0\equiv  1  {\ ({\rm mod}\ {6})}$, $U_l+T_l\equiv -1  {\ ({\rm mod}\ {6})}$ for $l\geq 2$.

This problem is connected with the primitive Pythagorean triples, balancing numbers and some more triangular numbers. Each number in the sequence (\ref{sequence1}) is the third (greatest) element $q$ of a Pythagorean triple of the form $(w,w+1,q)$, where $w\in\mathbb{N}_0$. Therefore, by (\ref{fraction}), we conclude that $1+24C$ is a sum of squares of the two consecutive natural numbers, i.e., \begin{equation}\label{kvocijent-w}C=\frac{w(w+1)}{12},\end{equation} which implies that \begin{equation}\label{&C-1}6C=\Delta_w,\end{equation}a triangular number. From (\ref{S'_11}) and (\ref{T'_11}), using (\ref{kvocijent-w}), it follows that $T'=2w+1$, $(S')^2=2w(w+1)+1$, and \begin{equation}\label{ppt}(w,w+1,q)=\Big(\frac{T'-1}{2},\frac{T'+1}{2},S'\Big)\end{equation} is a primitive Pythagoren triple. Since it is well known that the product of elements of every primitive Pythagorean triple is divisible by $60$, by using (\ref{T'_11}), we conclude that $5\mid C\cdot S'$. From (\ref{S'_11}), it follows that exactly one of $C$ and $S'$ is divisible by $5$. Let us also notice that the triple (\ref{ppt}) belongs to almost-isosceles right-angled triangle. By \cite[Theorem 2.1]{nyb}, all almost-isosceles right-angled triangles are given by \begin{equation}\label{ppt1}\Big(\frac{4D_n -1+\sqrt{8D_n^2+1}}{2}, \frac{4D_n +1+\sqrt{8D_n^2+1}}{2},2D_n +\sqrt{8D_n^2+1}\Big),\end{equation} where $D_n$ denotes the positive square root of the $n$-th square triangular number, so called balancing number (see \cite{panda}). By using (\ref{S'_11}), from (\ref{ppt}), (\ref{ppt1}), and \cite[Theorem 2.1]{panda}, we obtain \begin{equation}\label{quadratic}6C=D_n(3D_n+\sqrt{8D_n^2+1})=D_nD_{n+1}.\end{equation}From (\ref{quadratic}), (\ref{ppt}) and (\ref{ppt1}), we also have
\begin{eqnarray*}
T'-S'&=&2D_n,\\
T'+S'&=&2D_{n+1},
\end{eqnarray*}
where $D_{n+1}=3D_n+\sqrt{8D_n^2+1}$. Let us also notice (see \cite[Theorem 5.1(c)]{panda}) that $$S'T'=D_{2n+1}$$ is a balancing number, i.e., $(S'T')^2$ is a square triangular number.

Since the square of an odd number has the form $24l+1$ or $72l+9$, for some integer $l$, knowing that $c',d'\neq 0$, from (\ref{Z'_1}), we conclude that
\begin{eqnarray*}  c'd' &\equiv& 0{\ ({\rm mod}\ {3})}\ \ \textup{or}\\ c'd' &\equiv& 1{\ ({\rm mod}\ {3})},\end{eqnarray*}respectively.

\section{Fundamental solutions of the system of generalised Pellian equations}

Triangular $D(3\lambda)$-triples of the form (\ref{deset}) are quite often. For example, if $\lambda=C$ we have Triangular $D(3C)$-triples of the form
$\{3C,6C,24C+1\}$. One example of such a triple is a Triangular $D(105)$-triple
$\{105,210,841\}$, obtained for $k=7$, $\varepsilon =1$, $C=35$. We are interested in the egzistence of Triangular $D(3\lambda)$-quadruple which extend such a triple, so we want to solve the system of generalized Pellian equations  (\ref{ZX'})-(\ref{ZY'}).  We follow the theory of Nagell \cite[Theorem 108.a]{nagel} and adapt the result of Dujella \cite[Lemma 1]{Dabs} for the Triangular $D(3\lambda)$-quadruple. First we have to find the fundamental solutions of the associated Pell equations \begin{eqnarray}\label{petnaest1}U^2-3CV^2&=&1,\\
\label{sesnaest1}U^2-6CW^2&=&1.\end{eqnarray}

The equations (\ref{S'_11}) and (\ref{T'_11}) hold and, as in \cite{FJS}, using Magma we found that $3C+1$, $24C+1$ and $48C+1$ cannot be squares simultaneously (except for $C=0$). Similarly, we tested each one of $6C+1$, $12C+1$, and $27C+1$, together with  (\ref{S'_11}) and (\ref{T'_11}). It follows that the fundamental solutions\footnote{If we formulate the system (\ref{ZX'})-(\ref{ZY'}) as \begin{eqnarray*}
3\lambda(Z')^2- c(X')^2&=3\lambda-c,\\
6\lambda(Z')^2-c(Y')^2&=6\lambda-c,
\end{eqnarray*}which may seem simpler and easier to solve, then we would have associated Pell equations \begin{eqnarray*}U^2-3C(1+24\lambda)V^2&=&1,\\
U^2-6C(1+24\lambda)W^2&=&1\end{eqnarray*}for which we would not be able to determine fundamental solutions in general case.} of the equations (\ref{petnaest1}) and (\ref{sesnaest1}) are $(T',4)$ and $(S',2)$, respectively. Therefore,  the fundamental solution of \eqref{petnaest1} is
\[
U+V\sqrt{3C}=T'+4\sqrt{3C}.
\]
Thus, infinitely many solutions in natural numbers of the equation (\ref{ZX'}) are given by \begin{eqnarray}\label{osamnaest}Z'\sqrt{3\lambda}+X\sqrt{c'}=
(Z'_0\sqrt{3\lambda}+X_0\sqrt{c'})(T'+4\sqrt{3C})^m,\ \ m\geq 0,
\end{eqnarray}where $(Z'_0,X_0)$ is fundamental solution of the equation (\ref{ZX'}). The fundamental solution of the Pell equation  (\ref{sesnaest1})
is
\[
U+W\sqrt{6C}=S'+2\sqrt{6C}.
\]
Therefore, infinitely many solutions in natural numbers of the equation (\ref{ZY'}) are given by \begin{eqnarray}\label{devetnaest}Z'\sqrt{6\lambda}+Y\sqrt{c'}=
(Z'_1\sqrt{6\lambda}+Y_1\sqrt{c'})(S'+2\sqrt{6C})^n,\ \ n\geq 0,
\end{eqnarray}where $(Z'_1,Y_1)$ is fundamental solution of the equation (\ref{ZY'}).

From (\ref{osamnaest}) we conclude that $Z'=V_m$, for some non-negative integer $m$, where
\begin{eqnarray}\label{dvadesetjedan}V_0&=&Z'_0,\ \ V_1=T'Z'_0+4c'X_0,\ \ V_{m+2}=2T'V_{m+1}-V_m,\end{eqnarray}and from (\ref{devetnaest})  we conclude that $Z'=W_n$, for some non-negative integer $n$, where \begin{eqnarray}
\label{dvadesetdva}W_0&=&Z'_1,\ \ W_1=S'Z'_1+2c'Y_1,\ \ W_{n+2}=2S'W_{n+1}-W_n.\end{eqnarray}For any solution $(X,Y,Z')$ of the system (\ref{ZX'})-(\ref{ZY'}), we have \begin{equation}\label{dvadesetjedan1}Z'=V_m=W_n,\end{equation}
for some non-negative integers $m$ and $n$ and the sequences $(V_m)_{m\geq 0}$ and $(W_n)_{n\geq 0}$ given by (\ref{dvadesetjedan})-(\ref{dvadesetdva}), respectively. Equation (\ref{dvadesetjedan1}) leads to the extension of the Triangular $D(3\lambda)$-triple $\{3\lambda,6\lambda,c\}$ to a quadruple. In order to solve it, we will use some congruences obtained from recurrence sequences \eqref{dvadesetjedan} and \eqref{dvadesetdva}:
\begin{align}
V_{2m}&\equiv V_0=Z_0'\pmod{8c'},\quad W_{2n}\equiv W_0=Z_1'\pmod{4c'},\label{mod:even}\\
V_{2m+1}&\equiv V_1\pmod{8c'},\qquad~\, W_{2n+1}\equiv W_1 \pmod{4c'}.\label{mod:odd}
\end{align}We will also try to explore some additional properties of $Z'_0$, $Z'_1$, $X_0$ and $Y_1$ and describe these values more precisely.

If $\underline{Z'=V_{2m}}$, from (\ref{mod:even}) and (\ref{dvadesetjedan1}), it follows that $Z'\equiv Z'_0  {\ ({\rm mod}\ {8c'})}$ so $(Z')^2\equiv (Z'_0)^2  {\ ({\rm mod}\ {8c'})}$. Using (\ref{Z'}), we obtain \begin{eqnarray}\label{Z_0-new}(Z'_0)^2\equiv 1  {\ ({\rm mod}\ {8c'})}.\end{eqnarray}

If $\underline{Z'=V_{2m+1}}$, from  (\ref{osamnaest}) using (\ref{T'_1}), we obtain \begin{equation}\label{new}(Z'\sqrt{3\lambda}+X\sqrt{c'})(T'-4\sqrt{3\lambda c'})=(Z'_0\sqrt{3\lambda}+X_0\sqrt{c'})(T'+4\sqrt{3\lambda c'})^{2m}.\end{equation}By comparing the coefficients next to $\sqrt{3\lambda}$ on both sides of the equation (\ref{new}), we get \begin{equation}\label{eq-V_2m}T'Z'-4c'X=V_{2m}.\end{equation} By squaring (\ref{eq-V_2m}) and using (\ref{mod:even}), (\ref{dvadesetjedan1}) and (\ref{T'_1}), $(Z')^2\equiv (Z'_0)^2  {\ ({\rm mod}\ {8c'})}$ arises and (\ref{Z_0-new}) holds again. Since $Z'=V_{2m}$ or $Z'=V_{2m+1}$, by \eqref{Z_0-new}, we see that $d_0:=((Z_0')^2-1)/(8c')$ is a non-negative integer. Therefore, \begin{equation}\label{Z_0^2-new}(Z'_0)^2=8c'd_0+1.\end{equation}

If $\underline{Z'=W_{2n}}$, from  (\ref{mod:even}), $W_{2n}=Z'_1+4c'l$, for some $l\in\mathbb{Z}$. By squaring that and using  (\ref{Z'}), we get $ (Z'_1)^2 \equiv 1{\ ({\rm mod}\ {8c'})}.$  Therefore, \begin{equation}\label{Z_1^2-new1}(Z'_1)^2=8c'd_1+1,\end{equation}for some non-negative integer $d_1$.

If $\underline{Z'=W_{2n+1}}$, from  (\ref{mod:odd}), (\ref{dvadesetdva}) and (\ref{S'_1}), $(Z')^2\equiv (Z'_1)^2  {\ ({\rm mod}\ {4c'})}$. By (\ref{Z'}), we find $(Z'_1)^2\equiv 1  {\ ({\rm mod}\ {4c'})}$. Therefore, in this case, there exists non-negative $d_1\in\mathbb{Z}$ such that \begin{equation}\label{Z_1^2-new}(Z'_1)^2=4c'd_1+1.\end{equation} \\

By (\ref{Z_0^2-new}),  $(Z'_0)^2\geq 1$. By the equation (\ref{ZX'}), it holds \begin{equation}\label{ZX'-factor}3\lambda((Z'_0)^2-1)=c'(X_0^2-1-24\lambda).\end{equation}By inserting (\ref{Z_0^2-new}) into (\ref{ZX'-factor}), we further get \begin{equation}\label{X'_0^2-new}X_0^2=24\lambda(d_0+1)+1.\end{equation}Therefore, $X_0^2\geq 1+24\lambda$ and $X_0\geq 4k+\varepsilon$. Also, $X_0=4k+\varepsilon$ if and only if $Z'_0=\pm 1$. In that case $d_0=0$.

For $Z'_1$, it holds (\ref{Z_1^2-new1}) if $Z'=W_{2n}$ or (\ref{Z_1^2-new}) if $Z'=W_{2n+1}$. Hence, $(Z'_1)^2\geq 1$. By the equation (\ref{ZY'}), it holds \begin{equation}\label{ZY'-factor}6\lambda((Z'_1)^2-1)=c'(Y_1^2-1-24\lambda).\end{equation}By inserting (\ref{Z_1^2-new1}) into (\ref{ZY'-factor}), we get \begin{equation}\label{Y'_1^2-new}Y_1^2=24\lambda(2d_1+1)+1\end{equation}for the case where $Z'=W_{2n}$ and, by inserting (\ref{Z_1^2-new}) into (\ref{ZY'-factor}), we get \begin{equation}\label{Y'_1^2-new1}Y_1^2=24\lambda(d_1+1)+1\end{equation}for the case where $Z'=W_{2n+1}$. In both cases,  $Y_1^2\geq 1+24\lambda$ and $Y_1\geq 4k+\varepsilon$. Also, $Y_1=4k+\varepsilon$ if and only if $Z'_1=\pm 1$.  In that case $d'_1=0$.

\bigskip

In order to determine fundamental solutions of the system (\ref{ZX'}) and (\ref{ZY'}), we need suitable upper bounds (expressed in $c'$) on these fundamental solutions. Since right hand sides of equations (\ref{ZX'}) and (\ref{ZY'}) are negative integers, we follow the arguments of \cite[Theorem 108.a]{nagel} and, the same way as Dujella in \cite[Lemma 1]{Dabs}, we prove
\begin{eqnarray}
0<&X_0^2&\leq \frac{24\lambda((24\lambda+1) c'-3\lambda)}{T'-1}<\frac{24C(24\lambda+1)}{T'-1},\label{in:X0}\\
&(Z'_0)^2&\leq\frac{(T'-1)((24\lambda+1) c'-3\lambda)}{6\lambda}<\frac{T'(24\lambda+1)C}{6\lambda^2},\label{in:Z0}\\
0<&Y_1^2&\leq \frac{12\lambda((24\lambda+1) c'-6\lambda)}{S'-1}<\frac{12C(24\lambda+1)}{S'-1},\label{in:Y1}\\
&(Z'_1)^2&\leq\frac{(S'-1)((24\lambda+1) c'-6\lambda)}{12\lambda}<\frac{S'(24\lambda+1)C}{12\lambda^2}.\label{in:Z1}
\end{eqnarray}

Upper bounds (\ref{in:X0})-(\ref{in:Z1}) are not useful enough in determining fundamental solutions of the equations  (\ref{ZX'})-(\ref{ZY'}). In order to get more suitable ones,  we assume that $c'> 48\lambda$ and the possibility $c'\leq 48\lambda$ will be observed separately.

Let us assume that $c'>48\lambda$. By (\ref{C_l}) and the fact that $\lambda\geq 5$, it follows that $C>1200$ so $C\geq 40426$. It is not hard to see that there is no any suitable value $\lambda$ which divides $C=40426$. Therefore, $$C\geq 1373295.$$ From that, it also follows $$c'\geq 8119.$$ We further have \begin{equation}\label{T'-bound_lambda}T'=\sqrt{48\lambda c'+1}>48.0001\lambda\end{equation} and \begin{equation}\label{S'-bound_lambda}S'=\sqrt{24\lambda c'+1}>33.9417\lambda.\end{equation}Now, from (\ref{in:X0}), (\ref{in:Y1}) and (\ref{in:Z1}), we have
\begin{eqnarray}0<&X_0&< 0.5063c'\leq 0.10126C,\label{in:X0-new}\\
0<&Y_1&< 0.5964 c'\leq 0.11928C,\label{in:Y1-new}\\
&|Z_1'|&<1.1966c'\leq 0.23932C.\label{in:Z1-new}\end{eqnarray} Also, the condition \begin{equation}|Z_0'|<2c'\leq 0.4C,\label{in:Z0-new}\end{equation} by \eqref{in:Z0}, is equivalent to
\begin{align*}
\big(48\lambda c'+2T'+2-4(24\lambda+1)^2\big)(c')^2+24\lambda(24\lambda+1)c'-36\lambda^2>0,
\end{align*}
which holds as far as $c'>48\lambda$.

\bigskip

%%%%%%%%%%
In the following two lemmas we analyze the possibility that $c'\leq 48\lambda$ for a certain upper bound on $k$.
\begin{lem}\label{triples-1000}
Let $3 \le k<1000$, with $4k+\varepsilon$ prime and $1 \leq c'\leq 48\lambda=16(2k^2+\varepsilon k)$.
All the triples $(k,c',\varepsilon)$, satisfying \eqref{S'} and \eqref{T'}, are
\begin{align}\label{small}
(k,c',\varepsilon)\in&\{(3,5,1),(3,7,-1),(3,170,1),(3,238,-1),\notag\\
&~(7,1,1),(7,34,1),(21,4785,-1),(45,34177,1)\}.
\end{align}\end{lem}
\begin{proof}
The set \eqref{small} is easily obtained by checking the conditions in the assertion.
\end{proof}
%%%%%%%%%%

\bigskip

\begin{lem}\label{quadruples-1000}
For $(k,c',\varepsilon)=(45,34177,1)$ the fundamental solutions of the system  (\ref{ZX'})-(\ref{ZY'}), for which extensions of the Triangular Diophantine $D(3\lambda)$-triple $\{3\lambda,6\lambda,(24\lambda+1)c'\}$ with $3\lambda=2k^2+\varepsilon k$ are possible, are given by
$$V_{2m}=W_{2n},\ \ Z_0=Z_1=\pm 1,\ \ X_0=Y_1=4k+\varepsilon.$$
For other $(k,c',\varepsilon)$ satisfying (\ref{small}) there does not exist Triangular $D(3\lambda)$-quadruple $\{3\lambda,6\lambda,(24\lambda+1)c',d\}$.
\end{lem}
\begin{proof}
Let us consider the system of Pellian equations (\ref{ZX'})-(\ref{ZY'}) for $(k,c',\varepsilon)$ satisfying \eqref{small}, i.e. the systems
\begin{eqnarray*}21(Z')^2-5X^2&=&-824,\\
42(Z')^2-5Y^2&=&-803;\end{eqnarray*}
\begin{eqnarray*}15(Z')^2-7X^2&=&-832,\\
30(Z')^2-7Y^2&=&-817;\end{eqnarray*}
\begin{eqnarray*}21(Z')^2-170X^2&=&-28709,\\
42(Z')^2-170Y^2&=&-28688;\end{eqnarray*}
\begin{eqnarray*}15(Z')^2-238X^2&=&-28783,\\
30(Z')^2-238Y^2&=&-28768;\end{eqnarray*}
\begin{eqnarray*}105(Z')^2-X^2&=&-736,\\
210(Z')^2-Y^2&=&-631;\end{eqnarray*}
\begin{eqnarray*}105(Z')^2-34X^2&=&-28489,\\
210(Z')^2-34Y^2&=&-28384;\end{eqnarray*}
\begin{eqnarray*}861(Z')^2-4785X^2&=&-32963004,\\
1722(Z')^2-4785Y^2&=&-32962143;\end{eqnarray*}
and
\begin{eqnarray}4095(Z')^2-34177X^2&=&-1119668602,\label{syst1}\\
8190(Z')^2-34177Y^2&=&- 1119664507;\label{syst2}\end{eqnarray}
respectively. For the first seven systems, we used an algorithm written by Sz.~Tengely \cite{AT} in Sage, analogously as it is explained in \cite[Lemma 2.1.]{FJS}. In case of each that system, we get $Z'=\pm 1$, which means that $d=0$.

For the system (\ref{syst1})-(\ref{syst2}), using Python we checked which $d_0$ and $d_1$ are possible such that the upper bounds (\ref{in:X0}) and (\ref{in:Y1}) hold and $24\lambda(d_0+1)+1$ and $24\lambda(d_1+1)+1$ are squares. We get that $d_0,d_1\in\{0,6\lambda\}=\{0,8190\}$.
Therefore, using (\ref{Z_0^2-new}) and (\ref{X'_0^2-new}), we get
\begin{itemize}
  \item $d_0=0$, $(Z'_0,X_0)=(\pm 1, 4k+\varepsilon)$ or
  \item $d_0=6\lambda$, $(Z'_0,X_0)=(\pm T', 12\lambda+1)$.
\end{itemize}Since $Z'=V_{2m}$ or $Z'=V_{2m+1}$, from (\ref{dvadesetjedan1}), (\ref{mod:even}) and (\ref{mod:odd}) it follows
\begin{eqnarray}
Z'&\equiv& \pm 1\pmod{4c'}\label{d_0-0-1}\ \ \textup{or}\\
Z'&\equiv&\pm T'\pmod{4c'},\label{d_0-0-2}
\end{eqnarray}respectively. Also, using (\ref{Z_1^2-new}) (which, algebraically, covers the possibility  (\ref{Z_1^2-new1})) and (\ref{Y'_1^2-new}), we get
\begin{itemize}
  \item $d_1=0$, $(Z'_1,Y_1)=(\pm 1, 4k+\varepsilon)$ or
  \item $d_1=6\lambda$, $(Z'_1,Y_1)=(\pm S', 12\lambda+1)$.
\end{itemize}Since $Z'=W_{2n}$ or $Z'=W_{2n+1}$, from (\ref{dvadesetjedan1}), (\ref{mod:even}) and (\ref{mod:odd}), for $d_1=0$ it follows
\begin{eqnarray}
Z'&\equiv& \pm 1\pmod{4c'}\label{d_1-0-1}\ \ \textup{or}\\
Z'&\equiv&\pm S'\pm 2c'\pmod{4c'}\label{d_1-0-2},
\end{eqnarray}respectively, and for $d_1=6\lambda$ it follows
\begin{eqnarray}
Z'&\equiv& \pm S'\pmod{4c'}\label{d_1-6l-1}\ \ \textup{or}\\
Z'&\equiv&\pm 1+2c'\pmod{4c'},\label{d_1-6l-2}
\end{eqnarray}respectively. In case of  the system (\ref{syst1})-(\ref{syst2}), $c'=34177$, $\lambda=1365$, $S'=33461$ and $T'=47321$, so it is easy to check that every combination of the congruences (\ref{d_1-6l-1}), (\ref{d_1-6l-2}) with (\ref{d_0-0-1}), (\ref{d_0-0-2}) leads to a contradiction. Also, if we assume that one of the congruences  (\ref{d_0-0-1}), (\ref{d_0-0-2}) and one of the congruences (\ref{d_1-0-1}), (\ref{d_1-0-2}) simultaneously hold, we easily obtain a contradiction for all possible combinations of the signs $\pm$, except for $Z'=V_{2m}=W_{2n}$, with $d_0=d_1=0$.
\end{proof}

Now we may assume that $c' > 48\lambda$, for $3 \le k <1000$, and try to determine all fundamental solutions $(Z_0',X_0)$ and $(Z_1',Y_1)$ of the system  (\ref{ZX'})-(\ref{ZY'}) under that assumptions. In order to solve the equation (\ref{dvadesetjedan1}), which leads to the extension of the Triangular $D(3\lambda)$-triple $\{3\lambda,6\lambda,c\}$ to a quadruple, we need the following two lemmas. Note that in the following lemma the assumption $c' > 48\lambda$ doesn't change the result, so we omit it. In that sense, the statement of the Lemma \ref{quadruples-1000} for $k\in S$ is also proved here.

%Also, if $c' > 3\lambda$, then we get
%\begin{align}\label{in:Z01}
%|Z_0'|<4c',\quad |Z_1'|<4c'.
%\end{align}
%Indeed, the condition $|Z_0'|<4c'$ is equivalent to
%\[
%(192\lambda c'-576\lambda^2+8T'-48\lambda+7)(c')^2+6\lambda(24\lambda+1)c'-9\lambda^2>0,
%\]
%which holds as far as $c'\ge 3\lambda$ (note that we then have $T'=\sqrt{48\lambda c'+1}>12\lambda$). The case where $c'=3\lambda$ is not possible, as %we already explained.
%The other inequality $|Z_1'|<4c'$ can be shown in a similar way.

%
\begin{lem}\label{lem:XY<}For $(k,\varepsilon)\in S,$ where $S$ is the set (\ref{thm:main-S}), it holds $X_0\in\{4k+\varepsilon,12\lambda+1\}$ or $X_0\geq 24\lambda+1$ and $Y_1\in\{4k+\varepsilon,12\lambda+1\}$ or $Y_1\geq 24\lambda+1$.\end{lem}
\begin{proof}For any pair $(k,\varepsilon)$, if $d_0=0$ then, by (\ref{X'_0^2-new}), $X_0=4k+\varepsilon$. Also, if $d_0=6\lambda$, we have $X_0=12\lambda+1$. For $d_0\neq 0$, by using a program in Sage, we checked the pairs $(k,\varepsilon)$, such that  $3 \le k <1000$, $\varepsilon\in\{-1,1\}$, and $4k+\varepsilon$ is a prime. We assumed that $X_0< 24\lambda+1$ i.e., by (\ref{X'_0^2-new}), $0<d_0< 24\lambda+1$. For every $\lambda=\frac{2k^2+\varepsilon k}{3}$ and $d_0$, satisfying these conditions, we checked if $24\lambda(d_0+1)+1$ is a square. If that condition was satisfied, then we checked the elliptic curve which appears from the fact that $$(24c'\lambda+1)(48c'\lambda+1)(8c'd_0+1)$$ should be a square, after a substitution analogue as in \cite{FJS}. For each elliptic curve $E$ defined this way, we searched for integral points on it and, consequently,  possible values $c'$, using the command {\tt E.integral\_points(mw\_basis$=$E.gens(proof$=$True))}.

The program showed that for $(k,\varepsilon)\in S$ there is no integral points on the elliptic curve. For other elliptic curves, we were unable to obtain a provable Mordell-Weil basis, so for these cases we resorted to heuristic methods by passing {\tt proof$=$False} to {\tt E.gens}. The frequency of these fallbacks increased with the value of $k$. Using this method, we performed a brute-force search up to $k = 249$, yielding no valid integral points and thereby further strengthening Conjecture \ref{conjecture}. We also note that the search becomes exceedingly slow for larger values of $k$.

For $d_0=3\lambda$ we got singular curves, so we had to exclude that possibility separately. In case of $d_0=3\lambda$, by (\ref{X'_0^2-new}) and (\ref{lambda}), $8(2k^2+\varepsilon k)^2+8(2k^2+\varepsilon k)+1$ should be a square, which is not possible for the pairs $(k, \varepsilon)$ from the set $S$.

We did completely analogous checking for $Y_1<24\lambda+1$. For $d_1=0$, $Y_1=4k+\varepsilon.$ There is also a possibility that $Y_1=12\lambda+1$, which appears for $d_1=3\lambda$ if $Z'=W_{2n}$ and for $d_1=6\lambda$ if $Z'=W_{2n+1}$.  For $d_1=6\lambda$ if $Z'=W_{2n}$ and for $d_1=12\lambda$ if $Z'=W_{2n+1}$, we have a singular curve. That cases are excluded since $32(2k^2+\varepsilon k)^2+8(2k^2+\varepsilon k)+1$ can't be a square for the pairs $(k, \varepsilon)$ from the set $S$. The pairs $(k,\varepsilon)\in S$ passed the checkings both for $X_0$ and for $Y_1$.

\end{proof}
\begin{lem}\label{lem:XY>}
If $X_0\ge 12\lambda+1$ and $Y_1 \ge 12\lambda+1$, then
\begin{align*}
0<4c'X_0-T'|Z_0'|<4c',\quad 0<2c'Y_1-S'|Z_1'|<2c'.
\end{align*}
\end{lem}
\begin{proof}
The inequality $4c'X_0-T'|Z_0'|<4c'$ is equivalent to
\[
16(c')^2(X_0^2-2X_0+1)<(48\lambda c'+1)(Z_0')^2,
\]
which is rewritten as $X_0 \ge 12\lambda+1$, by using equation \eqref{ZX'} with $(Z',X)=(Z_0',X_0)$.
The inequality $4c'X_0-T'|Z_0'|>0$ easily follows from the inequality (\ref{in:Z0-new}) and equation \eqref{ZX'}.
The other inequalities are shown in a similar fashion.
\end{proof}
%%%%%%%%%%

\bigskip

Now we are solving the equation (\ref{dvadesetjedan1}), depending on the parity of indices $m$ and $n$, for $(k,\varepsilon)\in S$. By Lemma \ref{lem:XY<}, for such pairs $(k,\varepsilon)$, we have $$X_0\in\{4k+\varepsilon,12\lambda+1\}\ \ \textup{or}\ \ X_0\geq 24\lambda+1 \ \ \textup{and}\ \ Y_1\in\{4k+\varepsilon,12\lambda+1\}\ \ \textup{or}\ \ Y_1\geq 24\lambda+1.$$ Each combination of the parity of indices $m$ and $n$ is observed separately.

\bigskip

\textbf{Case 1)} $Z'=V_{2m}=W_{2n}$. \\

It follows from \eqref{mod:even}, (\ref{in:Z0-new}), and \eqref{in:Z1-new} that $Z_0'=Z_1'$. It holds (\ref{Z_0^2-new}) and, by \eqref{in:Z0-new}, we see that \begin{equation}0\leq d_0:=((Z_0')^2-1)/(8c')<\frac{c'}{2}\label{d_0-bound}\end{equation} is a non-negative integer such that $d_0<c=(24\lambda+1)c'$.
Since
\begin{align*}
8\cdot 3\lambda d_0+(24\lambda+1)&=X_0^2,\\
16\cdot 3\lambda d_0+(24\lambda+1)&=Y_1^2,\\
8cd_0+(24\lambda+1)&=\left((4k+\varepsilon)Z_0'\right)^2,
\end{align*}
if $d_0>0$ then $\{3\lambda,6\lambda,c,d_0\}$ would be a Triangular $D(3\lambda)$-quadruple with $d_0<c$, which contradicts the assumption that $c$ is minimal unless the quadruple is improper, which may happen if $d_0\in\{3\lambda,6\lambda\}$. The case $d_0=3\lambda$ is excluded in Lemma \ref{lem:XY<}. If $d_0=6\lambda$, we have $X_0=12\lambda+1$, but also $d_1=6\lambda$, which is excluded in Lemma \ref{lem:XY<}.
The possibility $d_0=0$ yields $$Z_0'=Z_1'=\pm1\ \ \textup{and}\ \  X_0=Y_1=4k+\varepsilon.$$

\bigskip

\textbf{Case 2)} $Z'=V_{2m+1}=W_{2n}$. \\

By (\ref{mod:even}), (\ref{mod:odd}) and (\ref{dvadesetjedan}), we have
\begin{equation}\label{kongr1}T'Z'_0 \equiv Z'_1 {\ ({\rm mod}\ {4c'})}.\end{equation}If $Z'_0=\pm 1$ ($X_0=4k+\varepsilon=\sqrt{24\lambda+1}$), then  $\pm T' \equiv Z'_1 {\ ({\rm mod}\ {4c'})}$. From (\ref{T'_1}), since $c'>48\lambda$, we conclude that $(T')^2<(c')^2+1$ and then\begin{equation}\label{T'-bound}T'<1.001c'.\end{equation} Therefore, (\ref{in:Z0-new}) and (\ref{T'-bound}) imply $Z'_1=\pm T'$. In this case (\ref{Z_1^2-new1}) holds, i.e., $d_1=6\lambda$ which is not possible for pairs $(k, \varepsilon)$ from the set $S$, as we saw in Lemma \ref{lem:XY<}.

If $Z'_1=\pm 1$, from (\ref{kongr1}) it follows $T'Z'_0 \equiv \pm 1 {\ ({\rm mod}\ {4c'})}$. By multiplying that congruence with $T'$ and by using (\ref{T'_11}), (\ref{in:Z0-new}), and (\ref{T'-bound}), we get $Z'_0=\pm T'$. By (\ref{Z_0^2-new}), $d_0=6\lambda$ and then, by (\ref{X'_0^2-new}), $X_0=12\lambda+1$. Also, from (\ref{Z_1^2-new1}), $d_1= 0$ and then, from (\ref{Y'_1^2-new}), $Y_1=\sqrt{24\lambda +1}$. Therefore, we have
 $$Z_0'=\pm T',\ \ Z_1'=\pm1,\ \ X_0=12\lambda+1,\ \ \textup{and}\ \  Y_1=\sqrt{24\lambda +1}=4k+\varepsilon.$$

Assume now that $|Z'_0|>1$ and $|Z'_1|>1$. Assume also that $X_0\geq 24\lambda+1$.
Then, by Lemma \ref{lem:XY>}, \eqref{kongr1}, and (\ref{in:Z1-new}), we have the following two possibilities:
\begin{itemize}
\item[(i)] $4c'X_0-T'|Z_0'|=|Z_1'|$,
\item[(ii)] $T'|Z_0'|-4c'(X_0-1)=|Z_1'|$.
\end{itemize}

If (i) occurs, put $d_1=((Z_1')^2-1)/(8c')$, which is an integer by \eqref{Z_1^2-new1} and, by (\ref{in:Z1-new}), satisfies \begin{equation}\label{d_1-bound}0\le d_1<0.179c'.\end{equation} We have
\begin{align*}
24\lambda d_1+(24\lambda+1)&=(T'X_0-12\lambda|Z_0'|)^2,\\
48\lambda d_1+(24\lambda+1)&=Y_1^2,\\
8cd_1+(24\lambda+1)&=((4k+\varepsilon)Z_1')^2.
\end{align*}
The minimality of $c'$ implies that $d_1\in \{0,3\lambda\}$. The case with $d_1=6\lambda$ is already excluded in Lemma \ref{lem:XY<}. If $d_1=0$ then, by (\ref{Z_1^2-new1}), $Z'_1=\pm 1$, which is not possible. If $d_1=3\lambda$, then $8(2k^2+\varepsilon k)^2+8(2k^2+\varepsilon k)+1$ should be a square, which is not possible for $(k,\varepsilon)$ from $S$.

%However, we then have
%\[
%4c'X_0-T|Z_0'|=4(12\lambda+1)c'-(48\lambda c'+1)=4c'-1\ne 1=|Z_1'|,
%\]
%which contradicts (i).

If (ii) occurs, the following congruences are obtained from \eqref{dvadesetjedan} and \eqref{dvadesetdva}:
\begin{align*}
V_{2m+1}&\equiv (-1)^m\cdot 4c'X_0 \pmod{T'},\\
W_{2n}&\equiv \begin{cases}(-1)^{n/2}Z_1'\pmod{T'}&\text{if~\,$2\mid n$},\\
(-1)^{(n-1)/2}\cdot 4c'Y_1 \pmod{T'}&\text{if~\,$2\nmid n$}.
\end{cases}
\end{align*}

In case $2 \mid n$, we see from $Z_1'\equiv \pm 4c'X_0\pmod{T'}$ that
\[
4c'(X_0-1)\equiv \pm 4c'X_0 \pmod{T'},
\]
which implies that $2X_0\equiv 1 \pmod{T'}$. Since $X_0\geq 24\lambda+1$, from (\ref{in:X0}) we further get $X_0<\frac{24C}{T'-1}$. Multiplying that equation by $2$ and using (\ref{T'_11}), we further obtain \begin{equation}2X_0-1<T'.\label{X_0-T'}\end{equation} Therefore, $X_0=\frac{1}{2}$, which is not possible.

%we have
%\[
%2X_0=T'+1.
%\]
%By \eqref{ZX'} we get
%\[
%(Z_0')^2=4(c')^2+\frac{T'c'}{6\lambda}-2c'+\frac{c'}{12\lambda}+1.
%\]
%Since $T'=\sqrt{48\lambda c'+1}>48\lambda$ by $c'>48\lambda$, we obtain $(Z_0')^2>4(c')^2$, which contradicts \eqref{in:Z01<}.

In case $2 \nmid n$, we see from $(-1)^m\cdot 4c'X_0\equiv (-1)^{(n-1)/2}\cdot 4c'Y_1 \pmod{T'}$ that \begin{equation}\label{c-n}X_0 \equiv \pm Y_1 \pmod{T'}.\end{equation} From (\ref{X'_0^2-new}) and (\ref{Y'_1^2-new}), by (\ref{d_0-bound}) and (\ref{d_1-bound}), we have $X_0<0.5001T'+1$ and $Y_1<0.43T'$.  Hence, by (\ref{c-n}), $X_0=\pm Y_1$. Since $X_0=-Y_1$ is not possible, by (\ref{c-n})
\[
X_0=Y_1.
\]
Then, it follows from equations \eqref{ZX'} and \eqref{ZY'} that
\[
(Z_0')^2=2(Z_1')^2-1.
\]
Inserting these equalities into
\[
(48\lambda c'+1)(Z_0')^2=(Z_1')^2+16(c')^2(X_0-1)^2+8c'(X_0-1)|Z_1'|,
\]
followed by (ii), and using (\ref{in:Z1-new}), we see that
\begin{align*}
X_0=Y_1&=\frac{2\big((24\lambda+2)c'-3\lambda\big)}{4c'-|Z_1'|}+\frac{1-(Z_1')^2-8c'|Z_1'|}{8c'(4c'-|Z_1'|)}\\
&<24\lambda+2-\frac{3\lambda}{c'}+\frac{1}{32(c')^2}<24\lambda+2,
\end{align*}
which contradicts the assumption $X_0\geq 24\lambda+1$ (note that $X_0$ is odd).

We still have to consider a possibility that  $X_0=12\lambda+1$, i.e $d_0=6\lambda$. The case with $d_0=0$ leads to $Z'_0=\pm 1$ and is already examined. The case with $d_1=0$ is already examined. If $d_1=3\lambda$, equations (\ref{Z_0^2-new}) and (\ref{Z_1^2-new1}) give $Z'_0=\pm T'$ and $Z'_1=\pm S'$. This, together with the congruence (\ref{kongr1}) and (\ref{T'_1}), gives $$\pm 1 \equiv \pm S' {\ ({\rm mod}\ {4c'})}.$$ By (\ref{S'_1}), we deduce \begin{equation}S'<0.71c',\label{S'-c'-bound}\end{equation}so we have $S'=\pm 1$, which is not possible. The last possibility is $Z'_0=\pm T'$ and $Y_1\geq 24\lambda+1$. From  (\ref{kongr1}) and (\ref{T'_1}), it follows $$\pm 1 \equiv Z'_1 {\ ({\rm mod}\ {4c'})}.$$ By (\ref{in:Z1-new}), $Z'_1=\pm 1$, which is a contradiction.

%%%%%%%%%%

\bigskip

\textbf{Case 3)} $Z=V_{2m}=W_{2n+1}$.\\

%%%%%%%%%%
By \eqref{dvadesetdva}, \eqref{mod:even} and \eqref{mod:odd}, we have
\begin{align}\label{kongr2}
Z_0'\equiv S'Z_1'+2c'Y_1 \pmod{4c'}.
\end{align}

Assume that $Y_1 \ge 12\lambda+1$. Then, by Lemma \ref{lem:XY>} and \eqref{kongr2}, we have
\begin{align}\label{eq:case3}
2c'Y_1-S'|Z_1'|=|Z_0'|.
\end{align}
We know that $d_0$ is an integer satisfying (\ref{d_0-bound}).
We then get
\begin{align*}
24\lambda d_0+(24\lambda+1)&=X_0^2,\\
48\lambda d_0+(24\lambda+1)&=(S'Y_1-12\lambda|Z_1'|)^2,\\
8cd_0+(24\lambda+1)&=((4k+\varepsilon)Z_0')^2.
\end{align*}
Since we are assuming $c'$ is minimal, we obtain $d_0\in\{0,6\lambda\}.$ The case with $d_0=3\lambda$ is excluded in Lemma \ref{lem:XY<}. The case where $d_0=6\lambda$ is not possible, because it would lead to a square of the form  $32(2k^2+\varepsilon k)^2+8(2k^2+\varepsilon k)+1$, which is not possible for the pairs $(k,\varepsilon)\in S$. Assume that $d_0=0$ and $Z_0'=\pm1$. Then, by \eqref{eq:case3} we have $|Z_1'|\equiv -S' \pmod{2c'}$, which together with (\ref{in:Z0-new}) and \eqref{S'-c'-bound} yields $|Z_1'|=2c'-S'$.
Again \eqref{eq:case3} implies that
\begin{align*}
2c'Y_1=|Z_0'|+S'|Z_1'|=1+S'(2c'-S')=2c'(S'-12\lambda),
\end{align*}
which yields $Y_1=S'-12\lambda$.
However, we then have
\begin{align*}
6\lambda(Z_1')^2-c'Y_1^2=6\lambda(2c'-S')^2-c'(S'-12\lambda)^2=6\lambda-c',
\end{align*}
which contradicts \eqref{ZY'}.

We also have to check the possibility $Y_1=4k+\varepsilon$. Then we have $d_1=0$ and, by \eqref{Z_1^2-new}, $Z'_1=\pm 1$. Since $Y_1>1$ is odd, in \eqref{kongr2}, we have $$Z_0'\equiv \pm S'+ 2c' \pmod{4c'}.$$ Hence, by \eqref{in:Z0-new} and \eqref{S'-c'-bound}, $Z'_0=\pm S'\pm 2c'$. Since \eqref{in:Z0-new} holds, $Z'_0=\pm (2c'-S')$. If $d_0=0$, i.e. $Z'_0=\pm 1$, then we would have $\pm S'=2c'\pm 1$, which is in contradiction with \eqref{S'-c'-bound}. If $d_0=6\lambda$, i.e. $Z'_0=\pm T'$, then we have $\pm T'\pm S'=\pm 2c'$. That is also not possible, because of \eqref{S'-c'-bound} and \eqref{T'-bound}. Hence $X_0\geq 24\lambda+1$, i.e. $d_0\geq 24\lambda+1$, by \eqref{X'_0^2-new}. By inserting $Z'_0$ and $X_0$ into \eqref{ZX'}, we get $c'-S'+6\lambda=2d_0\geq 48\lambda+2$, i.e. $c'\geq S'+42\lambda$. Together with \eqref{S'-bound_lambda}, that gives $(c')^2-84\lambda c'+42^2\lambda^2-1 \ge1$, which yields $c'>87.9411\lambda$ and $S'>45.9411 \lambda$.
Therefore, $S'<0.5246c'$, $64.9705\lambda < T' < 0.7404c'$ and $X_0 < 0.3194c'$. By \eqref{X_0-T'}, it also holds $X_0<\frac{T'+1}{2}$. In this case, we have
 \begin{equation}\label{case3}Z_0'=\pm (2c'-S'),\ \ Z_1'=\pm1,\ \ X_0\geq 24\lambda+1,\ \ \textup{and}\ \  Y_1=4k+\varepsilon.\end{equation}  Here, $X_0$ is attached to $Z_0$ by the equation \eqref{ZX'}. We have already noticed that \begin{equation}\label{d_0-3}d_0 = \frac{c' - S'}{2} + 3\lambda.
\end{equation}Since $S'$ is odd, we conclude that $c'$ should also be an odd integer. By substituting \eqref{d_0-3} into \eqref{X'_0^2-new} and using \eqref{S'_1}, we obtain
\begin{equation}
    2X_0^2 = (S' - 12\lambda)^2 + 48\lambda + 1. \label{eq:masterm}
\end{equation}  In Theorem \ref{thm:main}, we have \begin{equation}\label{lambdas}\lambda\in\{5, 7, 15, 35, 57, 77, 117, 495\}.\end{equation}All values $\lambda$ from \eqref{lambdas} are odd and  $c'$ should be odd, so $C$ should also be odd. Possible values of $C$ are elements of the sequence $(C_l)_{l\geq 0}$, given in \eqref{sequence2}. The sequence $C_l \pmod {2}$ cycles through $\{0, 1, 1, 0\}$ with period $4$. Since we need $2\nmid C_l$, we need $l \equiv 1,2 \pmod 4$ i.e. \begin{equation}\label{condition-l}l \equiv 1,2,5,6,9,10\pmod {12}. \end{equation} Sequences $C_l \pmod {5}$ and $C_l \pmod {7}$ both cycles through $\{0, 1, 0\}$ with period $3$. Hence, if $5\mid\lambda$ or $7\mid\lambda$, it follows that $5\mid C_l$ or $7\mid C_l$, respectively, which together with \eqref{condition-l} gives \begin{equation}\label{condition-l-1}l \equiv 2,5,6,9 \pmod {12}.\end{equation}The sequence $C_l \pmod {3}$ cycles through $\{0, 1, 2, 2, 1, 0\}$ with period $6$. Hence, if $3\mid \lambda$, it follows that $3\mid C_l$ which is, because of \eqref{condition-l} true for \begin{equation}\label{condition-l-2}l \equiv 5,6 \pmod {12}.\end{equation}
\begin{itemize}
\item For $\lambda=5$, in \eqref{eq:masterm} we have $2X_0^2 = S'^2 - 120S' + 3841,$ so $2X_0^2 \equiv S'^2 - S' + 5 \pmod 7$.
Possible values of $S'$ are elements of the sequence $(S'_l)_{l\geq 0}$, given in \eqref{sequence1}. The sequence $S'_l \pmod 7$ cycles through $\{1, 5, 1\}$ with period 3.
\begin{itemize}
    \item If $S'_l \equiv 1 \pmod 7$, then $2X_0^2 \equiv 5 \pmod 7$ so $X_0^2 \equiv 6 \pmod 7$, which is not possible.
    \item If $S'_l \equiv 5 \pmod 7$, then $2X_0^2 \equiv 4 \pmod 7$ so $X_0^2 \equiv 2 \pmod 7$. That is possible for $l \equiv 1 \pmod 3$, which contradicts \eqref{condition-l-1}.
\end{itemize}
\item For $\lambda=495$,  in \eqref{eq:masterm} we have $2X_0^2 = S'^2 - 11880S' + 35307361,$ so $2X_0^2 \equiv S'^2 - S' + 5 \pmod 7$ and we follow the same argument as for $\lambda=5$.
\item For $\lambda=117$, in \eqref{eq:masterm} it holds $2X_0^2= S'^2 - 2808S'+1976833.$ Therefore, $2X_0^2 \equiv S'^2 - S' + 5 \pmod 7$ and we follow the same argument as for $\lambda=5$, except that we have a contradiction with \eqref{condition-l-2}.
\item If $\lambda=15$, then \eqref{condition-l-2} holds. By \eqref{eq:masterm}, we have $2X_0^2 = S'^2 - 360S' + 33121$ and $2X_0^2 \equiv S'^2 -8S'\pmod {11}$.  The sequence $S'_l \pmod{11}$ has period 12:
\begin{equation}\label{11}\{1, 5, 7, 4, 6, 10, 10, 6, 4, 7, 5, 1\} \end{equation} so, by \eqref{condition-l-2}, we have to check only the possibility $S'_l \equiv 10 \pmod {11}$. Then $2X_0^2 \equiv 20 \pmod {11}$ so $X_0^2 \equiv 10\pmod {11}$, which is not possible.
\item If $\lambda=77$, then $77\mid C_l$ and $C_l \pmod {77}$ cycles through $\{0, 1, 35, 35, 1, 0\}$ with period $6$, so $l \equiv 0,5 \pmod 6$ i.e $l \equiv 0,5,6,11 \pmod {12}$. Together with \eqref{condition-l-1}, we have $l \equiv 5,6 \pmod {12}$.  By \eqref{eq:masterm}, we have $2X_0^2 = S'^2 -1848 S' + 857473$ so $2X_0^2 \equiv S'^2 -3S'+3 \pmod {5}$.  Since the sequence $S'_l \pmod {5}$ cycles through \begin{equation}\label{5}\{1,0,4,4,0,1\}\end{equation} with period $6$, we  have to check only the possibility $S'_l \equiv 1 \pmod {5}$. Then $2X_0^2 \equiv 1 \pmod {5}$ so $X_0^2 \equiv 3 \pmod {5}$, which is not possible.
\item If $\lambda=57$, then \eqref{condition-l-2} holds.  By \eqref{eq:masterm}, we have $2X_0^2 = S'^2 -1368S' + 470593$ so $2X_0^2 \equiv S'^2 -3S'+3 \pmod {5}$. The contradiction follows analogously as for $\lambda=77$. 

\item For $\lambda=7$, from \eqref{eq:masterm} we have that $S'=S_l$, where $(S_l)_{l\geq0}$ is sequence given by
$$S_0=84\pm1,\,S_1=136\pm3,\, S_{l+2}=6S_{l+1}-S_l-336.$$
Then, we easily get

$$S_l=\frac{13\sqrt{2}\pm1}{2}(3+2\sqrt{2})^l-\frac{13\sqrt{2}\mp1}{2}(3-2\sqrt{2})^l+84.$$

On the other hand, from (\ref{TS'}) we have that $S'=s_{\nu}$, where the sequence $(s_{\nu})_{\nu\geq0}$ is given by

$$s_0=1,\,s_1=5,\,s_{\nu+2}=6s_{\nu+1}-s_{\nu},$$

i.e.

$$s_{\nu}=\frac{\sqrt{2}+1}{2\sqrt{2}}(3+2\sqrt{2})^{\nu}+\frac{\sqrt{2}-1}{2\sqrt{2}}(3-2\sqrt{2})^l.$$

Now we have to solve $S_l=s_{\nu}$. It is easy to check that if we put $\alpha=3+2\sqrt{2}$, for $l\geq1$ we have

$$\alpha^{l+1}<S_l<\alpha^{l+2},$$ and for $\nu\ge1$ 
$$\alpha^{\nu-1}<s_{\nu}<\alpha^{\nu}.$$
It yields $\nu=l+2$ in $S_l=s_{\nu}$. But it is also easy to check that for $l\geq1$ we have $S_l<s_{l+2}$. It gives us that there is no solution of $S_l=s_{\nu}$, thus contradiction in this case.

\item The case $\lambda=35$ is similar to previous one. Firstly from \eqref{eq:masterm} we have that $S'=S_l^{(1),(2)}$, where $(S_l^{(1),(2)})_{l\geq0}$ is given by
$$S_0^{(1)}=420\pm1,\,S_1^{(1)}=536\pm3,\, S_{l+2}^{(1)}=6S_{l+1}^{(1)}-S_l^{(1)}-1680,$$ which yields

$$S_l^{(1)}=\frac{29\sqrt{2}\pm1}{2}(3+2\sqrt{2})^l-\frac{29\sqrt{2}\mp1}{2}(3-2\sqrt{2})^l+420,$$

or it is given by $$S_0^{(2)}=461,\,S_1^{(2)}=707,\, S_{l+2}^{(2)}=6S_{l+1}^{(2)}-S_l^{(2)}-1680,$$ i.e.

$$S_l^{(2)}=\frac{41\sqrt{2}+41}{2}(3+2\sqrt{2})^l-\frac{41\sqrt{2}-41}{2}(3-2\sqrt{2})^l+420.$$

As before from (\ref{TS'}) we have that $S'=s_{\nu}$, where the sequence $(s_{\nu})_{\nu\geq0}$ is given as in the previous case. Similarly it is easy to check that if we put $\alpha=3+2\sqrt{2}$, for $l\geq2$ we have

$$\alpha^{l+1}<S_l^{(1)}<\alpha^{l+2},$$ and for $\nu\ge1$
$$\alpha^{\nu-1}<s_{\nu}<\alpha^{\nu}.$$ It again gives $\nu=l+2$ in $S_l^{(1)}=s_{\nu}$ but for $l\geq2$ we have $S_l^{(1)}<s_{l+2}$. It gives us a contradiction, because it is easy to see that we do not have an intersection with small indices.

Furthermore, for $\l\geq2$ we have $$\alpha^{l+2}<S_l^{(2)}<\alpha^{l+3}.$$ This time it gives us $\nu=l+3$ in $S_l^{(2)}=s_{\nu}$ but for $l\geq2$ we have $S_l^{(2)}<s_{l+3}$, a contradiction again.

\end{itemize}

Therefore, \eqref{case3} is not an option for fundamental solutions of the system \eqref{ZX'}-\eqref{ZY'}.

\bigskip

\textbf{Case 4)} $Z=V_{2m+1}=W_{2n+1}$. \\

By \eqref{dvadesetjedan}, \eqref{dvadesetdva} and \eqref{mod:odd}, we have
\begin{align}\label{kongr3}
T'Z_0'+4c'X_0\equiv S'Z_1'+2c'Y_1 \pmod{4c'}.
\end{align}

Assume that $X_0\ge 12\lambda+1$ and $Y_1 \ge 12\lambda+1$. Then, Lemma \ref{lem:XY>} with \eqref{kongr3} shows that either of the following holds:
\begin{itemize}
\item[(i)] $4c'X_0-T'|Z_0'|=2c'Y_1-S'|Z_1'|$, with $Z_0'Z_1'>0$;
\item[(ii)] $4c'X_0-T'|Z_0'|=S'|Z_1'|-2c'Y_1+4c'$, with $Z_0'Z_1'<0$.
\end{itemize}

Consider first the case (i).
Put $Z''=4c'X_0-T'|Z_0'|=2c'Y_1-S'|Z_1'|$ and $d_2=((Z'')^2-1)/(8c')$.
Then, we see from \eqref{Z_1^2-new}, (\ref{Y'_1^2-new}), (\ref{T'_1}) and Lemma \ref{lem:XY>}, that $d_2$ is an integer satisfying $0 \le d_2<\frac{c'}{2}$.
As in Case 2) (i) or Case 3), we deduce that
\begin{align}\label{48}
24\lambda d_2+(24\lambda+1)&=(T'X_0-12\lambda|Z_0'|)^2,\notag\\
48\lambda d_2+(24\lambda+1)&=(S'Y_1-12\lambda|Z_1|)^2,\\
8cd_2+(24\lambda+1)&=((4k+\varepsilon)Z'')^2.\notag
\end{align}
By the minimality of $c'$, it should hold $d_2\in\{0,3\lambda,6\lambda\}$. The cases $d_2=3\lambda$ and $d_2=6\lambda$ can be excluded in the same fashion as we did for $d_0$ and $d_1$. Let $d_2=0$ and
\begin{align}\label{2cS}
Z''=4c'X_0-T'|Z_0'|=2c'Y_1-S'|Z_1'|=1.
\end{align}
We also see from \eqref{48} with $d_2=0$ and $(4k+\varepsilon)^2=24\lambda+1$ that
\begin{align}\label{S12}
S'Y_1-12\lambda|Z_1'|=4k+\varepsilon.
\end{align}
Combining \eqref{2cS} and \eqref{S12}, and using (\ref{S'-bound_lambda}), we get
\[
Y_1=\big((S')^2-24\lambda c'\big)Y_1=(4k+\varepsilon)S'-12\lambda>S'.
\]
But, by (\ref{Y'_1^2-new1}) and (\ref{d_0-bound}), $Y_1<S'$ and we have a contradiction.

Consider second the case (ii).
Note that \eqref{dvadesetjedan} and \eqref{dvadesetdva} imply that
\begin{align*}
V_{2m+1}&\equiv (-1)^m\cdot 4c'X_0 \pmod{T'},\\
W_{2n+1}&\equiv (-1)^{n(n+1)/2}\big(S'Z_1'+(-1)^n\cdot 2c' Y_1\big)\pmod{T'},
\end{align*}
whence, by $V_{2m+1}=W_{2n+1}$, we get
\begin{align*}
4c'X_0 \equiv \pm(S'Z_1'+(-1)^n\cdot 2c'Y_1)\pmod{T'}.
\end{align*}
It follows from (ii) that
\begin{align}\label{modT'}
\tau(S'Z_1'+(-1)^n\cdot 2c'Y_1)\equiv S'|Z_1'|-2c'Y_1+4c' \pmod{T'},
\end{align}
with $\tau \in \{-1,1\}$.

In case $\tau=-1$ and $2 \mid n$, congruence \eqref{modT'} shows that if $|Z_1'|=-Z_1'$, then $4c'\equiv 0 \pmod{T'}$, a contradiction; if $|Z_1'|=Z_1'$, then $S'|Z_1'|\equiv -2c' \pmod{T'}$, and by (ii),
\[
4c'X_0\equiv -2c'-2c'Y_1+4c' \pmod{T'},
\]
which yields $2X_0\equiv 1-Y_1 \pmod{T'}$.
By the inequalities from Lemma \ref{lem:XY>}, $2 \le 2X_0 \le T'+2$ and $Y_1<0.5983S'+1$, with $T'Y_1$ odd, we obtain $2X_0=2T'+1-Y_1$.
It follows from \eqref{ZX'} that
\begin{align*}
12\lambda (Z_0')^2&=\big(4X_0^2-(24\lambda+1)\big)c'+12\lambda\\
&=\big(16\cdot 12\lambda c'-(4T'+1-Y_1)(Y_1-1)-24\lambda+3\big)c'+12\lambda.
\end{align*}
We thus deduce from $Y_1<T'$ that
\begin{align*}
(Z_0')^2&=16(c')^2-\frac{(4T'+1-Y_1)(Y_1-1)c'}{12\lambda}-\frac{(8\lambda-1)c'}{4\lambda}+1\\
&>16(c')^2-\frac{(144\lambda c'+2-2T')c'}{12\lambda}-\frac{(8\lambda-1)c'}{4\lambda}+1\\
&=4(c')^2+\frac{(2T'+1-24\lambda)c'}{12\lambda}+1
>4(c')^2,
\end{align*}
which contradicts (\ref{in:Z0-new}).

In case $\tau=-1$ and $2 \nmid n$, congruence \eqref{modT'} shows that if $|Z_1'|=-Z_1'$, then $Y_1\equiv 1 \pmod{T'}$, which implies $Y_1=1$ and $(Z_1')^2=1-4c'$, a contradiction;
if $|Z_1'|=Z_1'$, then $2c'(Y_1-1)\equiv S'|Z_1'| \pmod{T'}$.
Combining this congruence with (ii), we get $4c'X_0 \equiv 2c'(Y_1-1)-2c'Y_1+4c'\pmod{T'}$, i.e., $2X_0 \equiv 1 \pmod{T'}$. As seen in (ii) of Case 2), if $X_0\geq 24\lambda+1$, then this leads to $X_0=\frac{1}{2}$, a contradiction. If $X_0=12\lambda+1$, by \eqref{T'-bound_lambda}, $2X_0-1<T'$, so we get $X_0=\frac{1}{2}$, again.

The cases $\tau=1$ can be shown in the same ways as the cases $\tau=-1$.

We are left with the possibilities that $X_0=4k+\varepsilon$ or $Y_1=4k+\varepsilon$, i.e. $d_0=0$ or $d_1=0$, respectively. The cases with $(d_0,d_1)\in\{(0,0),(0,6\lambda),(6\lambda,0)\}$ together with \eqref{kongr3} lead to $T'\equiv S' \pmod{2c'} $, $\pm 1\equiv S'\pmod{2c'}$ and $T'\equiv \pm1\pmod{2c'}$, respectively, which is a contradiction in all cases. If $d_0=0$ and $Y_1\geq 24\lambda+1$, by Lemma \ref{lem:XY>}, in \eqref{kongr3} we have $\pm T'\equiv S'Z'_1+2c'Y_1 \pmod{4c'}$ and then $\pm T'= 2c'Y_1- S'|Z'_1|.$ Let us denote $$Z''=\pm T'= 2c'Y_1- S'|Z'_1|.$$ By \eqref{T'_1}, we define $d'_0=\frac{(Z'')^2-1}{8c'}=6\lambda.$ Then\begin{align*}
24\lambda d'_0+(24\lambda+1)&=(12\lambda+1)^2,\\
48\lambda d'_0+(24\lambda+1)&=(S'Y_1-12\lambda|Z_1'|)^2,\\
8cd'_0+(24\lambda+1)&=((4k+\varepsilon)Z'')^2.
\end{align*}But, $d'_0=6\lambda$ is not possible for pairs $(k,\varepsilon)$ from the set $S$, as we already saw in similar situations.

If $d_1=0$ and $X_0\geq 24\lambda+1$, by Lemma \ref{lem:XY>}, in \eqref{kongr3} we have $T'Z'_0+4c'X_0\equiv  2c'\pm S' \pmod{4c'}$ and then $$2c'- S'=4c'X_0-T'|Z'_0|$$or $$4c'X_0-T'|Z'_0|-S'+ 2c'= 4c'.$$ Hence,
 \begin{equation}\label{case4}Z_1'=\pm 1,\ \ X_0\geq 24\lambda+1,\ \  Y_1=4k+\varepsilon,\ \ 4c'X_0-T'|Z'_0|=2c'\pm S',\end{equation} and equation \eqref{ZX'} holds. In order to exclude this possibility, we consider two sign-variants:

\medskip

\noindent\textbf{I.} $0<Z'_0<2c'$ and
\begin{equation}
4c'X_0 - T'Z'_0 = 2c'+\delta S' \quad\text{or}
\label{I}
\end{equation}

\noindent\textbf{II.} $-2c'<Z_0<0$ and
\begin{equation}
4c'X_0 + T'Z'_0 = 2c'+\delta S',
\label{II}
\end{equation}where \(\delta\in\{+1,-1\}.\) Note that variant II. is obtained from variant I. by the substitution \(Z'_0\mapsto -Z'_0\). Hence, it suffices to exclude variant I.

We already know that the admissible interval for $d_0$ is
\begin{equation}\label{range-d0}
24\lambda+1 \le d_0 < \frac{c'}{2}.
\end{equation}
Also, by (\ref{in:X0-new}),  \(X_0\) is an integer with \(24\lambda+1\le X_0\le \lfloor 0.5063c'\rfloor\). Assume one of the linear relations in \eqref{I} holds. Setting $R := 2c' + \delta S'$ and squaring, yields
\begin{equation}
16c'^2X_0^2 - 8c'X_0T'Z'_0 + T'^2(Z'_0)^2 = R^2. \label{0.7}
\end{equation} From \((\ref{X'_0^2-new})\),
\[
16c'^2X_0^2
= 384c'^2\lambda(d_0 + 1) + 16c'^2.
\] From \((\ref{T'_1})\) and \((\ref{Z_0^2-new})\),
\[
T'^2(Z'_0)^2
= 384c'^2\lambda d_0 + 48c'\lambda + 8c'd_0 + 1.
\]Finally, using (\ref{S'_1}), it follows
\[
R^2 = 4c'^2 + 4\delta c'S' + 24c'\lambda + 1.
\] Substituting these equalities into (\ref{0.7}) and dividing by \(8c'\), yields
\begin{equation}
X_0T'Z'_0
= 96c'\lambda d_0 + 48c'\lambda + \frac{3}{2}c'
+ 3\lambda + d_0 - \frac{\delta S'}{2}. \label{0.9}
\end{equation} On the other hand, from \((\ref{I})\) we have \(T'Z'_0 = 4c'X_0 - R\), so
multiplying by \(X_0\) and using (\ref{X'_0^2-new}) gives
\begin{equation}
X_0T'Z'_0= 96c'\lambda d_0 + 96c'\lambda + 4c' - RX_0. \label{0.10}
\end{equation}Equating (\ref{0.9}) and (\ref{0.10}), and substituting back \(R = 2c' + \delta S'\), we obtain
\begin{equation}
d_0 = 48c'\lambda + \frac{5c'+\delta S'}{2} - 3\lambda - (2c'+\delta S')X_0. \label{F}
\end{equation}
Therefore, $(5c'+\delta S')/2\in\mathbb Z$. Since $S'$ is odd, this forces $c'$ to be odd.

Let us consider any fixed pair $(\lambda,c')$, where $\lambda$ is from the set \eqref{lambdas} and $c'$ is any suitable value for that $\lambda$.  Combining (\ref{F}) with \eqref{range-d0} and solving for $X_0$ yields the real interval
\[
(X_0)_{\mathrm{low}} < X_0 \le (X_0)_{\mathrm{high}},
\]
where
\[
(X_0)_{\mathrm{high}} = \frac{48c'\lambda + \frac{5c'+\delta S'}{2} - 27\lambda-1}{2c'+\delta S'},
\qquad
(X_0)_{\mathrm{low}} = \frac{48c'\lambda + \frac{4c'+\delta S'}{2} - 3\lambda }{2c'+\delta S'}.
\] Since $c'>48\lambda+2$, by \eqref{range-d0}, and \(2c'+\delta S'> 1.29c'>0\), by \eqref{S'-c'-bound}, we obtain
\[
0< (X_0)_{\mathrm{high}}-(X_0)_{\mathrm{low}} =\frac{\frac{c'}{2} - 24\lambda-1}{2c'+\delta S'}< \frac{c'/2}{1.29c'} <0.388.
\]
Hence, the interval \(\langle (X_0)_{\mathrm{low}},(X_0)_{\mathrm{high}}]\) has length \(<1\) and contains at most one integer. Let us notice that it holds \begin{equation}\label{high}(X_0)_{\mathrm{high}} =24\lambda-\frac{24\lambda\delta S'+27\lambda+1}{2c'+\delta S'}+\frac{\frac{3}{2}c'}{2c'+\delta S'}+\frac{1}{2}.\end{equation}If  $\delta=1$, using $2c'+S'>2c'$, from \eqref{high} we obtain $$X_0<24\lambda+\frac{5}{4}.$$ Since $X_0\geq 24\lambda+1$, we conclude that $X_0=24\lambda+1$. By \eqref{X'_0^2-new},  $d_0=24\lambda+1$. If we now insert $(d_0,X_0)=(24\lambda+1,24\lambda+1)$ into equation \eqref{F}, we obtain \begin{equation}\label{eq1}54\lambda+2=c'+\delta S'(1+48\lambda)\end{equation}
which, together with  \eqref{S'_1}, form a system of two equations in unknowns $(c',S').$ Since in Theorem \ref{thm:main}, we have the set \eqref{lambdas} of values $\lambda$, we can solve the system of equations \eqref{eq1} and  \eqref{S'_1}, for each $\lambda$ from the set \eqref{lambdas} and neither in one case we get ineteger values for $c'$ and $S'$.

If  $\delta=-1$, from \eqref{high}, using \eqref{S'-c'-bound}, we obtain \begin{eqnarray}24\lambda+1<X_0 &<&24\lambda+\frac{24\lambda S'}{1.29c'}+\frac{\frac{3}{2}c'}{1.29c''}+\frac{1}{2}\nonumber\\ &<&37.20931\lambda+1.6628.\label{delta-1}\end{eqnarray}For each $\lambda$ from the set \eqref{lambdas}, inequalities \eqref{delta-1} give a range of possible values $X_0$. By using Sage, we checked which integers $d_0$ exist such that \eqref{X'_0^2-new} holds for $X_0$ from the range \eqref{delta-1}. For each such triple $(\lambda,d_0,X_0)$ the program solved the system of equations \eqref{S'_1} and \eqref{F} in $(c', S')$, and checked if $c'$ and $S'$ are integers. Neither one such pair of integers $(c',S')$ was found.

Hence, there is no admissible integer \(X_0\) for Variant I. Since we did not use the fact that $Z_0>0$ in the above argument, by the substitution \(Z_0\mapsto -Z_0\) Variant II is excluded as well, so the possibility \eqref{case4} for fundamental solutions is not possible.

\bigskip

We, therefore, have proved:
\begin{lem}\label{fundamentalna}Let $(k,\varepsilon)\in S$, where $S$ is the set from Theorem \ref{thm:main}. The fundamental solutions of the equations  (\ref{ZX'}) and (\ref{ZY'}), for which extensions of the Triangular Diophantine $D(3\lambda)$-triple $\{3\lambda,6\lambda,c\}$ are possible, are given by
\begin{itemize}
\item[\textbf{1.}] $V_{2m}=W_{2n}$, $Z'_0=Z'_1=\pm 1, X_0=Y_1=4k+\varepsilon,$
\item[\textbf{2.}] $V_{2m+1}=W_{2n}$, $X_0=12\lambda+1, Y_1=4k+\varepsilon, Z'_0=\pm T',  Z'_1=\pm 1,$

\end{itemize}
\end{lem}

\section{Linear forms in logarithms}

Now we are solving the equation $Z'=V_m=W_n$, where $m$ and $n$ are non-negative integers, for $(k,\varepsilon)\in S$, with the set $S$ from Theorem \ref{thm:main} and fundamental solutions of the equations  (\ref{ZX'}) and (\ref{ZY'}) described in Lemma \ref{fundamentalna}.  Using the standard method, we transform this equation to an inequality that includes a linear form in logarithms. We have $$V_m=\frac{1}{2\sqrt{3\lambda}}((Z'_0\sqrt{3\lambda}+X_0\sqrt{c'})(T'+4\sqrt{3\lambda c'})^m+(Z'_0\sqrt{3\lambda}-X_0\sqrt{c'})(T'-4\sqrt{3\lambda c'})^m)$$ and $$W_n=\frac{1}{2\sqrt{6\lambda}}((Z'_1\sqrt{6\lambda}+Y_1\sqrt{c'})(S'+2\sqrt{6\lambda c'})^n+(Z'_1\sqrt{6\lambda}-Y_1\sqrt{c'})(S'-2\sqrt{6\lambda c'})^n).$$ Introducing $$P:=\bigg(3\lambda Z'_0+\frac{X_0\sqrt{3\lambda c'}}{3\lambda}\bigg)(T'+4\sqrt{3\lambda c'})^m\ \ \textup{and} \ \   Q:=\bigg(6\lambda Z'_1+\frac{Y_1\sqrt{6\lambda c'}}{6\lambda}\bigg)(S'+2\sqrt{6c'\lambda})^n,$$
then $Z'=V_m=W_n$ implies $$P-\bigg(\frac{(24\lambda+1)c'}{3\lambda}-1\bigg)P^{-1}=Q-\bigg(\frac{(24\lambda+1)c'}{6\lambda}-1\bigg)Q^{-1}.$$ For $(Z'_0,X_0)$ and $(Z'_1,Y_1)$, from Lemma \ref{fundamentalna}, $P,Q>1$ holds for $c'\geq 8119$ and \begin{eqnarray*}P-Q&=&\bigg(\frac{(24\lambda+1)c'}{3\lambda}-1\bigg)P^{-1}-\bigg(\frac{(24\lambda+1)c'}{6\lambda}-1\bigg)Q^{-1}\\
&>&\bigg(\frac{(24\lambda+1)c'}{3\lambda}-1\bigg)(P^{-1}-Q^{-1})\\
&=&\bigg(\frac{(24\lambda+1)c'}{3\lambda}-1\bigg)(Q-P)P^{-1}Q^{-1}.\end{eqnarray*}Since $\frac{(24\lambda+1)c'}{3\lambda}-1>0$ and $P^{-1}Q^{-1}>0$, it follows that $P>Q$. Moreover, $$P-Q<\bigg(\frac{(24\lambda+1)c'}{3\lambda}-1\bigg)P^{-1},$$ so \begin{equation}\label{P-Q}\frac{P-Q}{P}<\bigg(\frac{(24\lambda+1)c'}{3\lambda}-1\bigg)P^{-2}.\end{equation}Then, using $c'\geq 8119$, $c'>48\lambda$ and Lemma \ref{fundamentalna}, for $m\geq 1$ we conclude
\begin{eqnarray*}P\geq \bigg(\frac{(12\lambda+1)\sqrt{3\lambda c'}}{3\lambda}-T'\bigg)(T'+4\sqrt{3\lambda c'})&>&\bigg(48\lambda c'-(T')^2+\frac{\sqrt{c'}}{\sqrt{3\lambda}}\bigg)\cdot 8\sqrt{3\lambda c'}\\
&=&8c'-8\sqrt{3\lambda c'}>6c',\end{eqnarray*}so $$P^{-1}<\frac{1}{6c'}.$$ For $c'\geq 8119$ and $\lambda\geq 5$, we also have $$\bigg(\frac{(24\lambda+1)c'}{3\lambda}-1\bigg)P^{-2}<\frac{1}{36c'}\bigg(\frac{121c'}{15}-1\bigg)<0.001.$$ Since $P>Q$, using (\ref{P-Q}) we get \begin{eqnarray*}0<\log \frac{P}{Q}&=&-\log\Big(1-\frac{P-Q}{P}\Big)\\
&<&\Big(\frac{(24\lambda+1)c'}{3\lambda}-1\Big)P^{-2}+\Big(\frac{(24\lambda+1)c'}{3\lambda}-1\Big)^2P^{-4}\\
&=& \Big(1+\Big(\frac{(24\lambda+1)c'}{3\lambda}-1\Big)P^{-2}\Big)\Big(\frac{(24\lambda+1)c'}{3\lambda}-1\Big)P^{-2}\\
&<&1.001\Big(\frac{(24\lambda+1)c'}{3\lambda}-1\Big)\frac{(T'+4\sqrt{3\lambda c'})^{-2m}}{\big(\frac{(12\lambda+1)\sqrt{3\lambda c'}}{3\lambda}-T'\big)^2}\\
&<&(T'+4\sqrt{3\lambda c'})^{-2m+2}.\end{eqnarray*} If we put $$\Lambda=m\log (T'+4\sqrt{3\lambda c'})-n\log (S'+2\sqrt{6\lambda c'})+\log\frac{2(3\lambda Z'_0+X_0\sqrt{3\lambda c'})}{6\lambda Z'_1+Y_1\sqrt{6\lambda c'}},$$ then we have proved:
\begin{lem}\label{m_i_n}
	If $V_m=W_n$  for $m,n\neq 0$, then for $(k,\varepsilon)\in S$, with the set $S$ from Theorem \ref{thm:main}, it holds\begin{eqnarray*}0<\Lambda<(T'+4\sqrt{3\lambda c'})^{-2m+2}.\end{eqnarray*}
\end{lem}

\section{Relation between indices $m$ and $n$ and their lower bound}

We want to find the relation between indices $m$ and $n$ if $V_m=W_n$, for $m,n\neq 0$. By (\ref{dvadesetjedan}) and (\ref{dvadesetdva}), we have $$V_m\geq (2T'-1)^{m-1}V_1\ \ \textup{and}\ \  W_n<(2S')^{n-1}W_1.$$ By  Lemma \ref{fundamentalna}, there are two possibilities for the fundamental solutions of equations (\ref{ZX'}) and (\ref{ZY'}):
\begin{itemize}
\item[\textbf{1.}] $V_{2m}=W_{2n}$, $Z'_0=Z'_1=\pm 1, X_0=Y_1=4k+\varepsilon$\\
In this case, $$(2T'-1)^{2m-1}(\pm T'+4(4k+\varepsilon) c')<(2S')^{2n-1}(\pm S'+2(4k+\varepsilon) c'),$$ which implies $$m<n.$$
\item[\textbf{2.}] $V_{2m+1}=W_{2n}$, $X_0=12\lambda+1, Y_1=4k+\varepsilon, Z'_0=\pm T',  Z'_1=\pm 1$\\
Then,  $$(2T'-1)^{2m}(\pm (48\lambda c'+ 1)+4(12\lambda+1)c')<(2S')^{2n-1}(\pm S'+2(4k+\varepsilon)c')$$
furthermore implies $$(2T'-1)^{2m}<\frac{2(4k+\varepsilon)c'-S'}{4c'-1}(2S')^{2n-1}<\frac{(8k+2)c'}{4c'-1}(2S')^{2n-1}<2.17k(2S')^{2n-1}.$$ Obviously, $2m\neq2n-1$. If we check the last inequality for $2m\geq 2n$, using \eqref{S'_1} and \eqref{lambda}, we find that it is not possible. Therefore, again we have $$m<n.$$
\end{itemize}

Now we intend to determine an upper bound on $n$ in terms of $m$. In order to do that we consider the "worst" possible case of fundamental solutions $(Z'_0,X_0)$ and $(Z'_1,Y_1)$ in the inequality of Lemma \ref{m_i_n}. Using the assumptions $\lambda\geq 5$ and $c'\geq 8119$, we obtain
\begin{eqnarray*}n&<&m\frac{\log (T'+4\sqrt{3\lambda c'})}{\log (S'+2\sqrt{6\lambda c'})}
+
\frac{
\log\!\left(
\frac{6\lambda T'+2(12\lambda+1)\sqrt{3\lambda c'}}
     {6\lambda+(4k+\varepsilon)\sqrt{6\lambda c'}}
\right)
}
{\log (S'+2\sqrt{6\lambda c'})}\\
&<&1.05m+0.361<1.42m.\end{eqnarray*}
Thus, we proved:

\begin{lem}\label{m<n}
	If $V_m=W_n$ for some $m,n\neq 0$, then  for $(k,\varepsilon)\in S$, with the set $S$ from Theorem \ref{thm:main}, it holds $m<n<1.42m.$
\end{lem}

In the sequel, we need the following congruences obtained by induction based on  (\ref{dvadesetjedan}) and (\ref{dvadesetdva}):
\begin{eqnarray*}V_{2m}&\equiv&(96m^2\lambda c'+1)Z'_0+8mT'c'X_0{\ ({\rm mod}\ {192\lambda c'^2})},\\
V_{2m+1}&\equiv&(96\lambda m(m+1)c'+1)T'Z'_0+4(2m+1)c'X_0{\ ({\rm mod}\ {192\lambda c'^2})},\\
W_{2n}&\equiv&(48n^2\lambda c'+1)Z'_1+4nS'c'Y_1{\ ({\rm mod}\ {192\lambda c'^2})},\\
W_{2n+1}&\equiv&(48\lambda n(n+1)c'+1)S'Z'_1+2(2n+1)c'Y_1{\ ({\rm mod}\ {192\lambda c'^2})}.\end{eqnarray*}
Again, we are separately considering the possibilities for the fundamental solutions of the system (\ref{ZX'}) - (\ref{ZY'}), given in Lemma \ref{fundamentalna}:
\begin{itemize}
\item[\textbf{1.}] $V_{2m}=W_{2n}$, $Z'_0=Z'_1=\pm 1, X_0=Y_1=4k+\varepsilon$\\

We have $$\pm(96m^2\lambda c'+1)+8(4k+\varepsilon)mT'c'\equiv\pm(48n^2\lambda c'+1)+4(4k+\varepsilon)nS'c'{\ ({\rm mod}\ {192\lambda c'^2})}.$$
If we divide it by $4c'$, we get \begin{equation}\label{congr}\pm24\lambda m^2+2(4k+\varepsilon)mT'\equiv\pm12\lambda n^2+(4k+\varepsilon)nS'{\ ({\rm mod}\ {48\lambda c'})}.\end{equation}
If $m\leq \frac{1}{2}\sqrt{c'}$, then congruence \eqref{congr} leads to a contradiction. Using the bounds $\frac{5}{3}k^2\leq 3\lambda$, $\lambda\geq 5$ and $c'\geq 8119$, we get
\begin{eqnarray*}24\lambda m^2&\leq&6\lambda c',\\
2(4k+\varepsilon)mT'&\leq&(4k+\varepsilon)\cdot\sqrt{c'}\cdot\sqrt{48\lambda c'+1}<5.814\sqrt{\lambda c'}\sqrt{48\lambda c'+1}<40.281\lambda c',\\
12\lambda n^2&<&12\lambda\cdot(1.42m)^2<6.05\lambda c',\\
(4k+\varepsilon)nS'&<&(4k+\varepsilon)1.42\cdot\frac{1}{2}\sqrt{c'}\sqrt{24\lambda c'+1}<20.222\lambda c'.\end{eqnarray*} Both sides of congruence (\ref{congr}) are positive because
\begin{eqnarray*}24\lambda m^2&\leq&24\lambda m\cdot\frac{1}{2}\sqrt{c'}=12\lambda m\sqrt{c'}<2(4k+\varepsilon)mT',\\
12\lambda n^2&<&12\lambda n\cdot1.42m< 17.05\lambda n\cdot \frac{1}{2}\sqrt{c'}<(4k+\varepsilon)nS'.\end{eqnarray*}
 Thus, we have equation $$\pm24\lambda m^2+2(4k+\varepsilon)mT'=\pm12\lambda n^2+(4k+\varepsilon)nS'.$$ With the $+$ sign, left-hand side is larger than the right-hand side. If the sign is $-$, then we also have the sign $-$  on right and \begin{eqnarray*}2(4k+\varepsilon)mT'-24\lambda m^2&=&m(2(4k+\varepsilon)T'-24m)>67.89\lambda m\sqrt{c'},\\
(4k+\varepsilon)nS'-12\lambda n^2&<&(4k+\varepsilon)nS'<(4k+\varepsilon)\cdot 1.42m \sqrt{24\lambda c'+1}<40.45\lambda m\sqrt{c'}.\end{eqnarray*} Hence, we obtained a contradiction and we may take $m>\frac{1}{2}\sqrt{c'}$ i.e. $$2m>\sqrt{c'}.$$
\item[\textbf{2.}] $V_{2m+1}=W_{2n}$, $X_0=12\lambda+1, Y_1=4k+\varepsilon, Z'_0=\pm T',  Z'_1=\pm 1$\\

Then  $$\pm(96\lambda m(m+1)c'+1)T'^2+4(12\lambda+1)(2m+1)c'\equiv\pm(48n^2\lambda c'+1)+4(4k+\varepsilon)nS'c'{\ ({\rm mod}\ {192\lambda c'^2})}$$
holds. Using the fact that in this case $Z_0Z_1>0$, we obtain $$\pm(96\lambda m(m+1)c'+48\lambda c')+4(12\lambda+1)(2m+1)c'\equiv\pm 48n^2\lambda c'+4(4k+\varepsilon)nS'c'{\ ({\rm mod}\ {192\lambda c'^2})}.$$ Dividing it by $4c'$,  we further get
$$\pm(24\lambda m(m+1)+12\lambda)+(12\lambda+1)(2m+1)\equiv\pm 12n^2\lambda+(4k+\varepsilon)nS'{\ ({\rm mod}\ {48\lambda c'})}.$$ Let us assume that $n\leq 0.71\sqrt{c'}$. Then $2m+1<2n\leq 1.42\sqrt{c'}$ so $m<0.71\sqrt{c'}$. Observe that
\begin{eqnarray*}24\lambda m(m+1)+12\lambda&<&12.289\lambda c',\\
(12\lambda+1)(2m+1)&<&0.042\lambda c',\\
12n^2\lambda&<&1.01\lambda c',\\
(4k+\varepsilon)nS'&<&20.222\lambda c'.\end{eqnarray*} So, the last congruence becomes equation
\begin{equation}\label{jednakost2}\pm(24\lambda m(m+1)+12\lambda)+(12\lambda+1)(2m+1)=\pm 12n^2\lambda+(4k+\varepsilon)nS'.\end{equation}With the sign $+$ on both sides of (\ref{jednakost2}), we have $24\lambda n(n+1)+12\lambda+(12\lambda+1)(2n+1)> 12n^2\lambda+(4k+\varepsilon)nS'$ i.e. $12\lambda n^2+(48\lambda+2)n+24\lambda+1>(4k+\varepsilon)nS'>24.09n\lambda \sqrt{c'}$, where the last inequality is obtained using $k>1.341\sqrt{\lambda}$. On the other hand, $12\lambda n^2+(48\lambda+2)n+24\lambda+1\leq 12\lambda n^2+(72\lambda+3)n\leq 12\lambda n^2+72.6\lambda n<9.33n\lambda \sqrt{c'}$. Therefore, we have a contradiction. If we have the sign $-$, then $-24\lambda m^2+2m+1=-12\lambda n^2+(4k+\varepsilon)nS'.$ We have $12\lambda n^2<12\lambda n\cdot0.71\sqrt{c'}\leq 8.52\sqrt{\lambda}n \sqrt{\lambda c'}<8.52\sqrt{\frac{2k^2+\varepsilon k}{3}}\cdot n\cdot\frac{1}{\sqrt{24}}S'<1.74knS'<(4k+\varepsilon)nS'$ and the left-hand side of \eqref{jednakost2} is negative for $m\geq 1$, so we have a contradiction. This implies $n> 0.71\sqrt{c'}$. Then, $$1.42(2m+1)>2n>2\cdot 0.71\sqrt{c'}$$ so $2m+1>\sqrt{c'}$.
\end{itemize}

We have, therefore, proved:
\begin{lem}\label{m_i_C}
	If $Z'=V_m=W_n$, for $m,n\neq 0$, then $m>\sqrt{c'}$.
\end{lem}

\section{Proof of the Theorem \ref{thm:main}}

To complete the proof of the Theorem \ref{thm:main}, we combine the following theorem with the result of Lemma \ref{m_i_n}. This theorem is due to Baker and W\"{u}stholz \cite{B-W}.
\begin{thm}(see \cite[Theorem]{B-W})\label{B-W-th}
	Let $\Lambda\neq 0$ be linear form in logarithms of $\ell$ algebraic numbers $\alpha_1,...,\alpha_{\ell}$ with rational integer coefficients $b_1,...,b_{\ell}$. Then $$\log \Lambda\geq -18(\ell +1)!\ell^{\ell+1}(32d)^{\ell+2}h'(\alpha_1)\cdot\cdot\cdot h'(\alpha_\ell)\log(2\ell d)\log B,$$ where $B=\max\{\mid b_1\mid,...,\mid b_\ell\mid\}$, and $d$ is the degree of the number field generated by  $\alpha_1,...,\alpha_{\ell}$. Here $$h'(\alpha)=\max\Big\{h(\alpha),\frac{\mid\log\alpha\mid}{d_\alpha},\frac{1}{d_\alpha}\Big\},$$ and $h(\alpha)$ denotes the standard Weil logarithmic height of $\alpha$, where $d_\alpha$ is the degree of the algebraic number $\alpha$.
\end{thm}

We take $$\alpha_1=T'+4\sqrt{3\lambda c'},\ \ \alpha_2=S'+2\sqrt{6\lambda c'},\ \ \alpha_3=\frac{Z_0'+\frac{X_0\sqrt{3\lambda c'}}{3\lambda}}{Z_1'+\frac{Y_1\sqrt{6\lambda c'}}{6\lambda}},$$ $\ell=3$, $B=n$, and $d=4$ in order to apply Theorem \ref{B-W-th} for the linear form $\Lambda$ which appears in Lemma \ref{m_i_n}. Because we are solving $Z'=V_m=W_n$ for $m,n\neq0$, as before, during the following calculations, we will use $\lambda\geq5$, $c'\geq8119$, $3\lambda=2k^2+\varepsilon k$ and $c'>48\lambda$.

Using the well-known methods, we get
$$h'(\alpha_1)=\frac{\log(T'+4\sqrt{3\lambda c'})}{2}=0.5\log\alpha_1,$$
$$h'(\alpha_2)=\frac{\log(S'+2\sqrt{6\lambda c'})}{2}<0.52\log c',$$
$$h'(\alpha_3)\leq \frac{1}{4}\log\bigg(36\lambda^2\Big(\frac{(4k+\varepsilon)^2 c'}{6\lambda}-1\Big)^2\cdot\frac{(T'+\frac{(12\lambda+1)\sqrt{3\lambda c'}}{3\lambda})^2}{\frac{(4k+\varepsilon)^2 c'}{6\lambda}-1}\bigg)<\log c'.$$
Form Lemma \ref{m_i_n} we have
\begin{equation}\label{Lambda-donja}-\log\Lambda>(2m-2)\cdot\log \alpha_1.\end{equation}
By Theorem \ref{B-W-th} we furthermore obtain
\begin{equation}\label{Lambda-gornja}-\log \Lambda\leq 18\cdot 4!\cdot 3^{4}\cdot 128^{5}\cdot 0.5\cdot 0.52\cdot \log\alpha_1 \cdot \log^2 c'\cdot\log(24)\cdot\log n.\end{equation}
Putting together (\ref{Lambda-donja}) and (\ref{Lambda-gornja}), using Lemma \ref{m<n} and Lemma \ref{m_i_C} we conclude
$$2m-2<9.935\cdot10^{14}\cdot\log^2(m^2)\cdot\log(1.42m),$$
which yields $$1<m<2.05\cdot 10^{20}.$$ Therefore, $c'<m^2<4.203\cdot 10^{40}$.

% which means that $k\leq 25$. Thus we proved the following proposition.
%\begin{prop}Let $k>25$ be an integer. Then the Triangular Diophantine triple $\{3,6,c_k\}$ cannot be extended to a quadruple with the larger element.
%\end{prop}

%
To finish the proof of Theorem \ref{thm:main}, we have applied a version of Baker-Davenport reduction (see \cite[Lemma 5]{DP}) to the linear form from Lemma \ref{m_i_n}, which is already well-known method in proving similar results. We have implemented it in GP-Pari and in all cases, after only one step of reduction, we gained $m\leq8$ which is small enough to get a contradiction (and then no extension to a triple exist except trivial with $d=0$), because we have already proved that $m\geq\sqrt{c'}\geq\sqrt{8119}>90$ if $m,n\neq0$. It is easy to see that the cases where $m=0$ or $n=0$ are not possible, because then we would have $d=0$ or $d$ is not an integer.

\bigskip

\section*{Acknowledgements}
The second author is supported by JSPS KAKENHI Grant Number JP24K06654. The first and the third author are supported by Croatian Science Fund, grant number IP-2022-10-5008.  The third author is also supported by European Union – NextGenerationEU - Particular problems in algebraic structures in number theory and topology. In the preparation of the proof of parts 3. and 4. of Lemma 4.5, the authors utilized ChatGPT-5.2 Thinking to assist in eliminating the cases \eqref{case3} and \eqref{case4}. All resulting solutions were corrected and subsequently verified manually by the authors. The authors take full responsibility for the mathematical correctness of these results.

%%%%%%%%%%%%%%%%%%%%%%%%%%%%%%%%%%%%%%%%%%%%%%%%%%%%%%%%%%%%%%%%%%%%%%

%%%%%%%%%%%%%%%%%%%%%%%%%%%%%%%%%%%%%%%%%%%%%%%%

\end{document}